\documentclass[smallextended]{svjour3}                     
\smartqed  
\usepackage{graphicx}
\usepackage{fix-cm}

\usepackage{endnotes}
\let\footnote=\endnote
\let\enotesize=\normalsize
\def\notesname{Endnotes}%
\def\makeenmark{$^{\theenmark}$}
\def\enoteformat{\rightskip0pt\leftskip0pt\parindent=1.75em
  \leavevmode\llap{\theenmark.\enskip}}

\usepackage{amsmath,amssymb,algpseudocode, multirow, float,bm,amsmath,bbm, theorem,accents,tikz,pgfplots}
\usepackage[bookmarksopen=true, bookmarksnumbered=true,
pdfstartview=FitH, breaklinks=true, urlbordercolor={0 1 0}, citebordercolor={0 0 1}]{hyperref}
\journalname{Mathematical Programming}
\pgfplotsset{compat=newest}

\usepgfplotslibrary{dateplot}
\usetikzlibrary{calc}
\newcommand{\md}[1]{\mathbb{#1}}
\newcommand{\ma}[1]{\mathcal{#1}}

\newcommand{\thickhline}{%
	\noalign {\ifnum 0=`}\fi \hrule height 1pt
	\futurelet \reserved@a \@xhline
}
\DeclareMathOperator*{\argmin}{arg\,min} 

\renewcommand{\P}{\mathbb{P}}
\newcommand{\E}{\mathbb{E}}

\newtheorem{prop}{Proposition}

\newtheorem{assumption}{Assumption}
\newtheorem{lem}{Lemma}

\newcommand*{\QED}{\hfill\ensuremath{\square}}
\newcommand{\one}{\mathbbm{1}}

\usepackage[toc,page]{appendix}

\renewcommand{\Comment}[2][.55\linewidth]{%
	\leavevmode\hfill\makebox[#1][l]{\#~#2}}

\DeclareMathOperator\Unif{Unif}

\DeclareMathOperator\sign{Sign}

\DeclareMathOperator\NLP{NLP}

\usepackage{algorithm}
\usepackage{natbib}
 \bibpunct[, ]{(}{)}{,}{a}{}{,}%
 \def\bibhang{24pt}%
\begin{document}

\title{Stochastic Cutting Planes for Data-Driven Optimization}

\author{Dimitris Bertsimas        \and
        Michael Lingzhi Li 
}

\institute{D. Bertsimas \at
              Sloan School of Management and Operations Research Center \\
              Massachusetts Institute of Technology\\
              	Cambridge, MA 02139\\
              \email{dbertsim@mit.edu}    \\       
           \and
           M. L. Li \at
              Operations Research Center\\
              Massachusetts Institute of Technology\\
              	Cambridge, MA 02139\\
              \email{mlli@mit.edu}
}
\date{Received: 12/30/2020 / Accepted: date}

\maketitle
\begin{abstract}
We introduce a stochastic version of the cutting-plane method for a large class of data-driven Mixed-Integer Nonlinear Optimization  (MINLO) problems. We show that under very weak assumptions the stochastic algorithm is able to converge to an $\epsilon$-optimal solution with high probability. Numerical experiments on several problems show that stochastic cutting planes is able to deliver a multiple order-of-magnitude speedup compared to the standard cutting-plane method. We further experimentally explore the lower limits of sampling for stochastic cutting planes and show that for many problems, a sampling size of $O(\sqrt[3]{n})$ appears to be sufficient for high quality solutions. 
\end{abstract}



%

\section{Introduction}
The cutting plane method, or the outer approximation method, has been a popular method to solve Mixed Integer Non-Linear Optimization  (MINLO) problems since its original introduction in \cite{duran1986outer} and \cite{fletcher1994solving}. It tackles problems  of the form:
\begin{align*}
    &\min_{\bm{z} \in \ma{Z},\bm{\theta} \in \Theta} f(\bm{z},\bm{\theta})\\
    &\text{subject to } g_i(\bm{z},\bm{\theta})\leq 0, \;\; i \in [m] ,
\end{align*}
where $\bm{z}$ and $\bm{\theta}$ are integer and continuous variables constrained to be in sets $\ma{Z}$ and $\Theta$, respectively, and $f, g_i$ are convex. The method aims to generate the tangent plane of both the objective and the constraints at an initial feasible solution $(\bm{z}_0,\bm{\theta}_0)$, and solve the resultant mixed integer linear optimization (MILO) instead in order to find an optimal  solution $(\bm{z}_1,\bm{\theta}_1)$. Then, the process is repeated where the tangent plane around $(\bm{z}_1,\bm{\theta}_1)$ is added until the solution stops updating (with a tolerance of $\epsilon$).  \cite{fletcher1994solving} proved that this method does indeed converge to the optimal solution in finite steps, and illustrated its practical effectiveness. 

With the rise in availability of  data in recent years, the cutting plane method has been applied to many data-driven optimization problems, ranging from risk minimization \citep{franc2009optimized}, matrix completion \citep{bertsimas2018interpretable}, to multiple kernel learning and inference in graphical models \citep{franc2011cutting}. 

However, as the data sizes continue to grow, cutting plane approaches face challenges in scaling as generating the full tangent planes usually involves evaluating the objective function and its derivative across all samples. With sample sizes easily reaching scales of $10^6$ or higher, this presents a computational difficulty for utilizing cutting-plane methods in large-scale real world applications. 

In this paper, we  introduce a stochastic variant of the cutting plane algorithm that aims to greatly increase the scalability of the cutting plane problem by multiple orders of magnitude. Specifically, instead of aiming to generate the exact tangent plane in each iteration, we utilize an independent random subset of the data to generate an approximation of the tangent plane. A parallel analogy of this idea in continuous optimization is stochastic gradient descent, where only a subset of the samples is considered to calculate a gradient, resulting in large performance increases that enables the optimization of millions of variables (such as in a convolutional neural network). 

Despite the simplicity of the idea, the authors could not find any detailed work exploring the theoretical properties and computational performance of this approach. In this work, we characterize the general convex data-driven optimization problems that stochastic cutting planes can be applied to, and demonstrate that the proposed method is able to converge to a near optimal solution with an exponentially vanishing failure probability. We illustrate that on multiple distinct problems (sparse regression, support vector machines and stochastic knapsack problems), stochastic cutting planes  record an order-of-magnitude speedup compared to standard cutting planes and compare competitively with state-of-the art specialized algorithms.

\subsection{Structure}
In Section \ref{sec:probformulation}, we
    introduce the general class of convex data-driven optimization problems to which we apply our framework. 
In Section \ref{sec:theory}, we prove that the stochastic cutting plane algorithm converges to an optimal solution of the original problem with exponentially vanishing failure probability under very weak conditions.
   In Section \ref{sec:experiments}, we present computational results on a variety of problems to show that the algorithm considerably outperforms the standard cutting plane method, and extends scalability by multiple orders of magnitude. 

\subsection{Literature}
 \cite{duran1986outer} first introduced  the exact cutting-plane method  in the context of mixed integer non-linear optimization (MINLO) problems that are linear in the integer variables.
 \cite{fletcher1994solving}    extended the framework to the general class of MINLO problems and provided   theoretical guarantees for the case that the MINLO is convex. 
\cite{parikh1976approximate} investigated an approximate cutting-plane method that involved generating a sequence of hyperplanes using non-optimal dual solutions to the inner NLO problem. 
The approximate cutting planes generated are always valid for the original objective function, but  the approximation is only effective ``near an  optimal solution". 

Many other works entitled ``stochastic cutting planes" or similar are primarily focused on randomized constraint generation in the presence of a large (possibly infinite) number of constraints, such as \cite{volkov1997general} or \cite{calafiore2005uncertain} in robust optimization. Concretely, they consider problems  of the form:
\begin{align*}
    &\min_{\bm{z} \in \ma{Z},\bm{\theta} \in \Theta} f(\bm{z},\bm{\theta})\\
    &\text{subject to } g_i(\bm{z},\bm{\theta})\leq 0, \;\; i \in [m] ,
\end{align*}
where $m$ is very large or an infinite set. The proposed algorithms generally first solves the problem with only a random finite subset of constraints, and then add the remaining constraints iteratively (and randomly) till convergence. In particular, they assume that the function $f(\bm{z},\bm{\theta})$ is relatively simple to evaluate. This work focuses on the case where the objective function takes significant time to calculate and proves theoretical guarantees on optimal solution retrieval using stochastic cutting planes. Furthermore, we demonstrate real-world scalability with data beyond $10^6$ samples.  

In stochastic optimization, sample average approximation is frequently utilized to provide an approximation to the objective function (which is usually of the form of an expected value). For stochastic  MINLO problems, there has been work (see e.g. \cite{wei2004sample,kleywegt2002sample}) adapting the cutting planes methodology for the stochastic problem by solving the sample average approximation at every iteration using a new sample of size $N$ multiple times. This work complements previous work by exploring how to solve the sample average approximation with sample size $N$ effeiciently, when $N$ is large. In particular in Section \ref{sec:experiments}, we demonstrate how stochastic cutting planes generates an $100$x speedup to solve the sample average approximation version of the static stochastic knapsack problem (SSKP). 

 \cite{bertsimas2018interpretable}  utilized a similar idea to speed up the cutting-plane method utilized to solve a variant of the matrix completion problem. The theoretical guarantees in this paper were limited to the specific application and in addition required strong convexity, which is not required in the proof for the present paper.

\section{Problem Formulation and Examples}
\label{sec:probformulation}
We consider a general class of mixed-integer data-driven optimization problems with integer variables $\bm{z} \in \md{Z}^{p_1}$ and continuous variables $\bm{\theta} \in \md{R}^{p_2}$. In addition, there is data $\bm{d}_1,\ldots,\bm{d}_N$ obtained either through historical results or simulations. As common in these settings, the goal is to find a representation of the data under some model driven by the variables $\bm{z}, \bm{\theta}$.  

 Utilizing the shorthand $[N]$ to represent $\{1,\ldots,N\}$, we can write the objective function of such data driven problem as:
\[f(\bm{z},\bm{\theta};[N]):=f(\bm{z},\bm{\theta};\bm{d}_1,\ldots,\bm{d}_N),\]
where we use the notation $f(\bm{z},\bm{\theta};[N])$ to denote that the objective function is created with samples in the set $[N]$. One of the most common examples of this type of objective function is data-driven risk minimization, where we have $\bm{d}_i=(\bm{x}_i,\bm{y}_i)$, and the objective function is:
\[f(\bm{z},\bm{\theta};[N])= \frac{1}{N}\sum_{i=1}^N l(\bm{y}_i,m(\bm{x}_i, \bm{z},\bm{\theta})). \]
Here $l(\cdot,\cdot)$ is some pre-defined loss function describing the difference between the target $\bm{y}_i$ and the model fit $m(\bm{x}_i, \bm{z},\bm{\theta})$. 

We   consider the data-driven optimization problem:
\begin{align}
    &\min_{z \in \ma{Z}, \theta \in \Theta} f(\bm{z},\bm{\theta};[N])\label{eq:obj_func} \vspace{3pt}\\ 
    &\text{subject to } g_i(\bm{z},\bm{\theta})\leq 0,\;\;\; i\in [m], \nonumber 
\end{align}
with the following assumptions:
\begin{assumption}[Convexity]
\label{ass:convexity}
$f(\bm{z},\bm{\theta};[N])$ and $g_i(\bm{z},\bm{\theta})$ are at least once differentiable and convex in $\bm{z},\bm{\theta}$ for all $i \in [m]$. 
\end{assumption}
\begin{assumption}[Compactness]
\label{ass:compact}
$\ma{Z}$ is a finite set   in $\md{Z}^{p_1}$ and $\Theta$ is a compact set in $\md{R}^{p_2}$.
\end{assumption}
\begin{assumption}[Data Concentration]
\label{ass:random}
Let $S^n$ be a random subset of $[N]$ of size $n$ selected without replacement. Then with probability at least $1-\epsilon$, we have:
\begin{align*}
    |f(\bm{z},\bm{\theta};S^n)-f(\bm{z},\bm{\theta};[N])|&\leq M\sqrt{\frac{\log(\frac{1}{\epsilon})}{n}}, \vspace{3pt}\\
    \|\nabla f(\bm{z},\bm{\theta};S^n)-\nabla f(\bm{z},\bm{\theta};[N])\|&\leq M'\sqrt{\frac{(p_1+p_2)\log(\frac{1}{\epsilon})}{n}},
\end{align*}
where $M$ and $M'$ are absolute constants independent of $p_1,p_2, n, N$. 
\end{assumption}
The convexity and compactness assumptions are standard in the classic cutting-plane literature. In real-world problems, the convexity assumption might prove restrictive as many important problems do not admit convex formulations in both the integer and continuous variables. Nevertheless, the proposed algorithm can also be applied to problems that are non-convex in $f$ and $g$, but the solution would not carry any theoretical guarantees.

The data concentration assumption essentially states that the objective function (and its derivative) created with a random subset of the data is  concentrated around the value of the objective (and its derivative) under the full set of the data. In essence, this requires that a small instance of the problem is ``similar" to the original  problem. Note that this assumption in particular does \emph{not} require $f(\bm{z};\bm{\theta};S^n)$ to be an unbiased estimator of $f(\bm{z};\bm{\theta};N)$, as long the estimate $f(\bm{z};\bm{\theta};S^n)$ is sufficiently close.
 This is true for the vast majority of data-driven problems, and we now introduce a few prominent classes of these problems. 

\subsection{Data-Driven Risk Minimization}
\label{subsec:riskmin}
Given features $\bm{x}_i$ and responses $y_i$, $i \in [N]$, the data-driven risk minimization problem aims to solve the following unconstrained optimization problem:
\begin{equation}
\label{eq:central}
\argmin_{\bm{z} \in \ma{Z}, \bm{\theta} \in \Theta} \frac{1}{N} \sum_{i=1}^N l(\bm{y}_i, \bm{x}_i, \bm{z}, \bm{\theta}), 
\end{equation}
where $l$ is a function convex in $\bm{z}, \bm{\theta}$. This structure contains a vast array of problems, including:
\begin{itemize}
    \item \emph{Machine Learning}: This includes for example support vector machines: 
  $$l(\cdot,\cdot,\cdot,\cdot)=\max  \{1-y_i\bm{\theta}^T\bm{x}_i,0\},$$
  least-squares regression:
  $$l(\cdot,\cdot,\cdot,\cdot)=\|y_i-\bm{\theta}^T\bm{x}_i\|^2,$$ 
and  many others. 
    \item \emph{Treatment Effect Estimation}: This includes a wide range of data-driven optimization problems in pricing and policy treatment effect estimation (see \citealp [e.g.][]{rivers1988limited,blundell2004endogeneity,alley2019pricing}) where the loss is constructed from a parametric or semi-parametric loss function and there are further constraints on the parameters (non-negativity, bounds, etc). Many of these problems might involve non-convex loss functions, but they are usually solved by methods that assume convexity (e.g. gradient descent, cutting planes) or create successive convex approximations. For example, given a probit model $\hat{f}(\bm{x}_i,\bm{\theta})$ for estimating propensity of purchasing a model, the form of the loss function would appear as:
    \[l(\cdot,\cdot,\cdot,\cdot)=\frac{1}{N} \sum_{i=1}^N y_i\log(\hat{f}(\bm{x}_i,\bm{\theta}))+(1-y_i)\log(1-\hat{f}(\bm{x}_i,\bm{\theta})),\]
    where $N$ is the number of samples, and $(\bm{x}_i,y_i)$ are the historical data.
    \item \emph{Sample Average Approximation }: Sample Average Approximation (SAA) methods are commonly used in stochastic optimization problems where the master problem:
    \[\min_{\bm{\theta} \in \Theta}\E_{\bm{W}\sim P}[G(\bm{\theta},\bm{W})]\]
    is solved by  generating $(\bm{w}^i)_{i\in [N]}$, a sample of size $N$ (either taken from the distribution or created through historical samples), and solving the sample average approximation:
        \[\min_{\bm{\theta} \in \Theta}\frac{1}{N}\sum_{i=1}^NG(\bm{\theta},\bm{w}^i). \]
    The approximation is often solved many times with different samples to create confidence intervals. It is often assumed that the resulting problem is tractable, but it is not necessarily so when $N$ is large \citep{razaviyayn2014successive}. For a detailed discussion of the sample average approximation for stochastic approximation, please see e.g. \citet{kleywegt2002sample}.
\end{itemize}

We next show that Problem \eqref{eq:central}   satisfies the data concentration Assumption \ref{ass:random}.
\begin{lemma}
\label{lem:riskmin}
Let $l$ be a function that is at least once differentiable in $\bm{z}$ and $\bm{\theta}$. Then $f(\bm{z},\bm{\theta};[N])=\frac{1}{N}\sum_{i=1}^N l(\bm{y}_i, \bm{x}_i, \bm{z}, \bm{\theta}) $ satisfies Assumption \ref{ass:random}.  
\end{lemma}
\begin{remark}
There is no need to assume  convexity of the loss function $l(\cdot,\cdot)$ for it to satisfy the data concentration Assumption \ref{ass:random}.
Such a  convexity assumption is only necessary to prove global convergence results.
\end{remark}
\proof{Proof of Lemma \ref{lem:riskmin}:}
With $f(\bm{z},\bm{\theta};[N])=\frac{1}{N}\sum_{i=1}^N l(\bm{y}_i, \bm{x}_i, \bm{z}, \bm{\theta})$, we can invert Assumption \ref{ass:random} so that it states for all $\epsilon>0$, we need to prove that there exists universal constants $M,M'$ such that:
\begin{align*}
    \P\left(\left|\frac{1}{N}\sum_{i=1}^N l(\bm{y}_i, \bm{x}_i, \bm{z}, \bm{\theta})-\frac{1}{n}\sum_{i \in S^n} l(\bm{y}_i, \bm{x}_i, \bm{z}, \bm{\theta})\right|>\epsilon \right)&\leq \exp\left(-\frac{n\epsilon^2}{M}\right),  \vspace{3pt}\\
    \P\left(\left\|\frac{1}{N}\sum_{i=1}^N \nabla l(\bm{y}_i, \bm{x}_i, \bm{z}, \bm{\theta})-\frac{1}{n}\sum_{i \in S^n} \nabla l(\bm{y}_i, \bm{x}_i, \bm{z}, \bm{\theta})\right\|>\epsilon \right)&\leq \exp\left(-\frac{n\epsilon^2}{M'(p_1+p_2)}\right).
\end{align*}
These two statements follow directly from Hoeffding's inequality in sampling without replacement (see e.g. \cite{bardenet2015concentration}). 
\QED
\endproof

\subsection{Sparse Data-Driven Risk Minimization}
In  Section \ref{subsec:riskmin} we included many key problems in data-driven optimization, but the format required the loss function to be separable in the space of samples. 
In this section, we  give an example of a data-driven problem in which the objective function is not directly separable in the  space of samples, and show that Assumption \ref{ass:random} still holds.

In particular, we  consider the problem of  sparse linear  regression, where we seek  at most $k$ out of $p$ factors in the design matrix $\bm{X} \in \md{R}^{N \times p}$ to explain the response data $\bm{y} \in \md{R}^N$.  The optimization problem can be written as:
\begin{equation}
\label{eq:sparse_reg_orig}
    \min_{\bm{\beta}: \|\bm{\beta}\|_0\leq k} \;\;\frac{1}{N}\left(\|\bm{y}-\bm{X}\bm{\beta}\|_2^2 + \frac{1}{\gamma}\|\bm{\beta}\|_2^2\right).
\end{equation}

\cite{bertsimas2020sparse} showed that it can be reformulated as:
\begin{equation}
\label{eq:sparse_reg}
    \min_{\bm{z} \in \ma{Z}} \;\;\frac{1}{N}\bm{y}^T \left(I_N + \gamma \sum_{i=1}^pz_i \bm{X}_i\bm{X}_i^T\right)^{-1}\bm{y},
\end{equation}
where here $\ma{Z}=\{z \in \{0,1\}^p \mid \sum_{i=1}^p z_i=k\}$, $\bm{X}_i$ is the $i$th column of the design matrix, and $\gamma$ controls the regularization. Note that this optimization function satisfies Assumptions \ref{ass:convexity} and \ref{ass:compact} as the objective function is convex in $\bm{z}$, and the set $\ma{Z}$ is indeed finite.  
\begin{lemma}
\label{lem:sparsereg}
The function 
$f(\bm{z},\bm{\theta};[N])=\frac{1}{N}\bm{y}^T \left(I_N + \gamma \sum_{i=1}^pz_i \bm{X}_i\bm{X}_i^T\right)^{-1} \bm{y}$ satisfies Assumption \ref{ass:random}.  
\end{lemma}

The proof of the lemma is contained in Appendix \ref{app:lem_sparsereg_proof}. The key characteristic that allowed this problem to also satisfy Assumption \ref{ass:random} is that Assumption \ref{ass:random} does not require the random sample estimate to be an \emph{unbiased} estimate of the true loss function (and in this case it is indeed not true). \cite{bertsimas2019unified} further discusses how to transform other sparse problems to a convex integer optimization problem  that could be treated similarly. \cite{bertsimas2018interpretable} discusses a similar problem in matrix completion where the data could be sampled in two dimensions.

\section{Stochastic Cutting Planes and Theoretical Results}
\label{sec:theory}
In this section, we present a formal description of the cutting plane algorithm, its stochastic cutting plane variant and show the global convergence of the
stochastic cutting plane algorithm. The cutting plane algorithm  introduced in \cite{duran1986outer} and generalized in \cite{fletcher1994solving} 
    for Problem \eqref{eq:obj_func} iteratively rotates between solving a nonlinear subproblem and a master mixed-integer linear program (MILP) containing linear approximations of the functions involved. We first define the NLP subproblem as followed:
    \begin{align*}
        \NLP(\bm{z},[S_n]):=&\min_{\bm{\theta}\in \Theta} f(\bm{z},\bm{\theta};[S_n])\\&\text{subject to } g_i(\bm{z},\bm{\theta})\leq 0, \;\; i \in [m] 
    \end{align*}
    Then we present the  cutting plane algorithm in Algorithm \ref{alg:cutting_plane_full}. At every iteration, we first check if the current solution $\bm{z}_t, \bm{\theta}_t,\eta_t$ is still a feasible solution for the $t$th iteration, and return $\bm{z}_t, \bm{\theta}_t,\eta_t$ (as the optimal solution) if it is feasible. As the cutting plane iteratively adds constraints, if the previous solution is still feasible for the next iteration, then it is also optimal for the next iteration. 
    
    We further note that this guarantees that when the algorithm terminates at iteration $t=T$, the final cutting plane is at the optimal solution $(\bm{z}^*, \bm{\theta}^*)$ of the $T$th iteration (i.e. $\bm{z}_T=\bm{z}^*,\bm{\theta}_T=\bm{\theta}^*$), which would simplify the proof of some of our technical results. This is also how modern solvers implement cutting planes through lazy constraints \cite{gurobi}. 

\begin{algorithm}[ht] 
	\begin{algorithmic}[1]
		\Procedure{CUTPLANES}{$\{\bm{d}_i\}_{i=1}^N$}
		\State $t \gets 1$
		\State $\bm{z}_1,\bm{\theta}_1 \gets \text{Initialization}$ \Comment{Heuristic Warm Start}
		\State $\eta_1 \gets lb$ \Comment{Initialize feasible solution variable with lower bound}
		\While{$\eta_t<f(\bm{z}_t,\bm{\theta}_t; [N])$ \textbf{or} $\exists j \in [m]\;0 < g_j(\bm{z}_t,\bm{\theta}_t)$}\\
		 \Comment{\parbox[t]{.5\linewidth}{while the current solution is not feasible for the $t$th iteration}}
		\State $\begin{aligned}\bm{z}_{t+1}, \eta_{t+1} \gets \displaystyle &\argmin_{\bm{z}\in \ma{Z}, \bm{\theta}\in \Theta, \eta\geq lb} \;\; \eta\\  \text{s.t.} \quad &\eta \geq f(\bm{z}_i,\bm{\theta}_i; [N])+\nabla f(\bm{z}_i,\bm{\theta}_i; [N])^T\binom{\bm{z}-\bm{z}_i}{\bm{\theta} - \bm{\theta}_i}, \quad \forall i \in [t]\\ \quad & 0 \geq g_j(\bm{z}_i,\bm{\theta}_i)+\nabla g_j(\bm{z}_i,\bm{\theta}_i)^T\binom{\bm{z}-\bm{z}_i}{\bm{\theta} - \bm{\theta}_i}, \quad \forall i \in [t], \; \; \forall j \in  [m] \end{aligned}$ \\\Comment{Solve the master problem}
		 \State $\bm{\theta}_{t+1} \gets \arg_{\bm{\theta}}  \;\NLP(\bm{z}_{t+1},[N])$\Comment{Solve the NLP subproblem}
		\State $t \gets t+1$
		\EndWhile
		\State $\bm{z}^*, \bm{\theta}^*, \eta^*\gets \bm{z}_t, \bm{\theta}_t, \eta_t$		
		\State \textbf{return} $\eta^*, \bm{z}^*, \bm{\theta}^*$\Comment{Return the filled matrix $\bm{X}$}
		\EndProcedure
	\end{algorithmic}
	\caption{Cutting-plane algorithm. }
	\label{alg:cutting_plane_full}
\end{algorithm}

In many examples such as risk minimization and sparse regression, the cutting plane algorithm (and its variants) has produced state-of-the-art results. However, one important drawback is its scalability with respect to data. It was noted in \cite{franc2008optimized} that the scalability of the cutting-plane method for Support Vector Machines (SVMs), although significantly better than other exact methods, trails significantly behind stochastic methods such as stochastic gradient descent. This behavior is also observed in matrix completion in \cite{bertsimas2018interpretable}, where the full cutting plane algorithm scales less favorably than competing algorithms. 

To understand the scalability of cutting plane algorithm, observe that at every iteration we are only solving a  MILO problem, which has enjoyed great advances in recent years due to various algorithmic improvements (see \cite{bixby2007progress, achterberg2013mixed} for details). As a result, the amount of time spent in the actual solver is usually not excessive. On the contrary, the cutting plane generation for $f$ and $g$ usually accounts for the majority of the running time. In data-driven optimization problems defined in Section \ref{sec:probformulation}
this is greatly exacerbated by the fact that the function $f(\bm{z},\bm{\theta};[N])$ and its derivative $\nabla f(\bm{z},\bm{\theta};[N])$ needs to be evaluated for all samples in $[N]$, which could easily exceed  $10^6$
samples  or higher. 

Therefore, with the objective  to improve scalability for cutting-plane methods, we propose that instead of evaluating the full function and its derivative, we   randomly select $n$ samples without replacement in a set $S^n$ to evaluate the function and its derivative at every iteration of the cutting plane. When $n \ll N$, we obtain  a considerable speedup for the generation of the cutting plane. Concretely, we present the stochastic cutting plane 
algorithm in Algorithm \ref{alg:cutting_plane_sto}.
\begin{algorithm}[ht] 
	\begin{algorithmic}[1]
		\Procedure{Stochastic Cutting Planes}{$\{\bm{d}_i\}_{i=1}^N$,$\epsilon$}
		\State $t \gets 1$
		\State $\bm{z}_1,\bm{\theta}_1 \gets \text{Initialization}$ \Comment{Heuristic Warm Start}
		\State $\eta_1 \gets lb$ \Comment{Initialize feasible solution variable with lower bound}
		\State $S^n_1 \gets \text{random $n$-sized subset of } \{1,\ldots,N\}$ \Comment{Initialize random subset of samples}
		\While{$\eta_t<f(\bm{z}_t,\bm{\theta}_t; [S^n_t])$ \textbf{or} $\exists j \in [m]\;0 < g_j(\bm{z}_t,\bm{\theta}_t)$}\\
		 \Comment{\parbox[t]{.5\linewidth}{while the current solution is not feasible for the $t$th iteration}}
		\State $\begin{aligned}\bm{z}_{t+1},  \eta_{t+1} \gets \displaystyle &\argmin_{\bm{z}\in \ma{Z}, \bm{\theta}\in \Theta, \eta\geq lb} \;\; \eta\\  \text{s.t.} \quad &\eta \geq f(\bm{z}_i,\bm{\theta}_i; S^n_i)+\nabla f(\bm{z}_i,\bm{\theta}_i; S^n_i)^T\binom{\bm{z}-\bm{z}_i}{\bm{\theta} - \bm{\theta}_i}, \quad \forall i \in [t]\\\quad & 0 \geq g_j(\bm{z}_i,\bm{\theta}_i)+\nabla g_j(\bm{z}_i,\bm{\theta}_i)^T\binom{\bm{z}-\bm{z}_i}{\bm{\theta} - \bm{\theta}_i}, \quad \forall i \in [t], \; \; \forall j \in  [m] \end{aligned}$\\\Comment{Solve the master problem}
		\State $S^n_{t+1} \gets \text{random $n$-sized subset of } \{1,\ldots,N\}$ 
		 \State $\bm{\theta}_{t+1} \gets \arg_{\bm{\theta}}  \;\NLP(\bm{z}_{t+1},[S^n_{t+1}])$\Comment{Solve the NLP subproblem}
		\State $t \gets t+1$
		\EndWhile
		\State $\bm{z}^*, \bm{\theta}^*, \eta^*\gets \bm{z}_t, \bm{\theta}_t, \eta_t$		
		\State \textbf{return} $\eta^*, \bm{z}^*, \bm{\theta}^*$\Comment{Return the filled matrix $\bm{X}$}
		\EndProcedure
	\end{algorithmic}
	\caption{Stochastic Cutting-plane algorithm. }
	\label{alg:cutting_plane_sto}
\end{algorithm}

\subsection{Theoretical Guarantees}
\cite{fletcher1994solving} proved that the cutting plane algorithm is able to converge to an optimal solution in finite steps. In this section, we  show that the stochastic cutting plane algorithm, similarly, is able to converge to an $\epsilon$-optimal solution in finite steps with high probability. This is presented in the following theorem:
\begin{theorem}
\label{theo:convergence}
Let $\bm{d}_1,\ldots,\bm{d}_n$ be given data, and assume that Problem  (\ref{eq:obj_func}) is feasible and has an  optimal solution  $\bm{z}^*,\bm{\theta}^*$. Under  Assumptions \ref{ass:convexity}-\ref{ass:random},  Algorithm \ref{alg:cutting_plane_sto} converges after a finite number of iterations  $K$ with probability one, and $\E[K]\leq (e-1)|\ma{Z}|+1$. When it terminates, it produces a feasible solution $\tilde{\bm{z}}^*,\tilde{\bm{\theta}}^*$ that satisfies:
\begin{equation}
    f(\tilde{\bm{z}}^*,\tilde{\bm{\theta}}^*;[N])\leq  f(\bm{z}^*,\bm{\theta}^*;[N])+ \epsilon
\end{equation}
with probability at least $1- (2K+1)\exp\left(\frac{-n\epsilon^2}{(2+\sqrt{p_1+p_2})^2J}\right)$, where $J$ is an absolute constant. 
\end{theorem}
The proof of Theorem \ref{theo:convergence} is contained in Appendix \ref{app:theo_converge_proof}. This theorem shows that the stochastic cutting plane algorithm finds  an $\epsilon$-optimal solution with failure probability that is exponentially decaying with respect to the number of samples involved. We can see that the failure probability also decays exponentially in the optimality gap $\epsilon$, which indicates that the solutions to the stochastic cutting plane method are concentrated very close to the true solution, as long as a reasonable $n$ is selected.

We further note that although in the worst case, the expected number of cutting planes needed could be very high ($(e-1)|\ma{Z}|+1$), this bound is only a small constant factor away from the worst case for the deterministic cutting plane algorithm ($|\ma{Z}$) as shown in \citep{duran1986outer}. Furthermore, although the number of cutting planes needed could be exponential in $p_1$ (if $\bm{z}$ are all binary, then $\ma{Z}=2^{p_1}$),  in practice it usually converges much quicker \citep{fletcher1994solving}. Furthermore, in most application settings, we have $p_1\ll n \ll N$, so even with the worst-case scenario of $K\sim \exp(p_1)$, the failure probability is still on the order of $\exp\left(p_1-O(n\epsilon^2)\right)\ll 1$.

\section{Computational Results}
\label{sec:experiments}
In this section, we present computational results from a variety of examples that demonstrate the versatility, scalability and accuracy of the stochastic cutting plane algorithm \ref{alg:cutting_plane_sto}.
\subsection{Sparse Regression}
\label{subsec:sparsereg}

In this section, we apply Algorithm \ref{alg:cutting_plane_sto} to the sparse linear regression problem \eqref{eq:sparse_reg}. 
As shown in Lemma \ref{lem:sparsereg}, the sparse linear regression problem satisfies the data concentration assumption and thus   the stochastic cutting plane algorithm \ref{alg:cutting_plane_sto} applies. 
We explore the scalability and accuracy of Algorithm \ref{alg:cutting_plane_sto}.

To study the effectiveness of Algorithm  \ref{alg:cutting_plane_sto} we generate data  according to 
$\bm{y}=\bm{X}\bm{\beta}+\bm{\epsilon}$ where the entries of the design matrix $X_{ij} \sim N(0,1)$ independently. A random set $S$ of size $k$ is chosen from $\{1,\ldots,p\}$, and the coefficients are chosen so that $\bm{\beta}_j = 0,\; \;\;\forall j \not \in S$, and $ \bm{\beta}_j  \sim N(0,1),\; \;\; \forall j \in S$ independently. This ensures the exact sparsity of rank $k$ for the resulting model. The noise is sampled independently with $\epsilon_i \sim N(0,\sigma)$. 

For each parameter combination $(N,p,k,\sigma)$, we   record the time $T$ and the Mean Average Percentage Error (MAPE) of both the standard cutting-plane algorithm (Algorithm \ref{alg:cutting_plane_full})
 and the stochastic cutting plane algorithm  (SCP, Algorithm \ref{alg:cutting_plane_sto}). 
The $\gamma$ parameter for both algorithms are selected through cross-validation on a validation set of size $N$ using the same setup as above. We set the tolerance in the cutting plane algorithm to $\epsilon=10^{-4}$. The testing set is of size $N$ sampled using the true $\beta$ and independent $X_{ij}$ and $\epsilon_i$ with the same setup as above. Table \ref{tab:sparsereg} shows the mean results over 10 random generations of the dataset of the size $(N,p,k,\sigma)$, separated into 4 blocks where each variable is varied over a range of possibilities. For these experiments, we used $n=\min\{N, 10\sqrt{N}\}$ as the sampling size for stochastic cutting planes. In Section \ref{subsec:samplesize}, we  detail further experiments on varying the sample size.  The testing environment has a six-core i7-5820k CPU with 24GB of RAM. Both algorithms are implemented in pure Julia with Julia v1.4.2.  

\begin{table}[ht!]
	\centering
	\resizebox{0.9\columnwidth}{!}{%
		\begin{tabular}{|c|l|l|l|l||l|c|l|c|}
			\hline \multirow{2}{*}{} & \multirow{2}{*}{$\bm{N}$}&\multirow{2}{*}{$\bm{p}$}&\multirow{2}{*}{$\bm{k}$} &\multirow{2}{*}{$\bm{\sigma}$}& \multicolumn{2}{c|}{\textbf{Cutting Planes}} & \multicolumn{2}{c|}{\textbf{SCP}} \\\cline{6-9}
			&  &   & & & $\bm{T}$ & \textbf{MAPE} & $\bm{T}$ & \textbf{MAPE}   \\\hline
			\multirow{4}{*}{$N$} & $10^3$ & $10^2$ & 10 & 0.1& 2.4s & $4.4\%$& 1.9s & $4.4\%$ \\
			&	$10^4$ & $10^2$ & 10 & 0.1 & 5.7s & $3.9\%$& 3.7s & $3.9\%$\\
				&	$10^5$ & $10^2$ & 10 & 0.1& 170s & $4.5\%$& 7.8s & $4.5\%$\\
				&	$10^6$ & $10^2$ & 10 & 0.1& 1356s & $3.1\%$& 46s & $3.1\%$\\\hline
			\multirow{3}{*}{$p$} &	$10^5$ & $10^2$ & 10 & 0.1& 170s & $4.5\%$& 7.8s & $4.5\%$\\
			&	$10^5$ & $10^3$ & 10 & 0.1& 2073s & $4.1\%$& 98.5s & $4.1\%$ \\
				&	$10^5$ & $10^4$ & 10 & 0.1& 29770s & $4.7\%$& 945s & $4.7\%$\\\hline
			\multirow{3}{*}{$k$} & $10^5$ & $10^2$ & 10 & 0.1& 170s & $4.5\%$& 7.8s & $4.5\%$\\
			&	$10^5$ & $10^2$ & 20 & 0.1& 1621s & $2.9\%$& 41.9s & $2.9\%$ \\
				&	$10^5$ & $10^2$ & 50 & 0.1& 30608s & $1.8\%$& 937s & $1.8\%$\\\hline
			\multirow{3}{*}{$\sigma$} & $10^5$ & $10^2$ & 10 & 0.1& 170s & $4.5\%$& 7.8s & $4.5\%$\\
			&	$10^5$ & $10^2$ & 10 & 0.2& 197s & $6.7\%$& 8.9s & $6.7\%$\\
				&	$10^5$ & $10^2$ & 10 & 0.3& 231s & $9.8\%$& 10.7s & $9.8\%$\\\hline\hline

		\end{tabular}
	}
	\caption{Comparison of cutting planes and stochastic cutting planes on synthetic data for sparse regression.  }
	\label{tab:sparsereg}
\end{table}
We see that as $N$ grows, the advantage of stochastic cutting planes in using $n=O(\sqrt{N})$ samples to compute the cutting plane quickly increases. At the largest scale with $N=10^6$, stochastic cutting planes generates over a $30\times$ speedup compared to the standard cutting plane algorithm. Furthermore, we see that the stochastic cutting plane converges to the exact same solution as the standard cutting planes, as reflected the identical MAPE figures. This is a reflection of the exponentially vanishing failure probability as indicated in Theorem \ref{theo:convergence}.  Such order of magnitude speedup is also preserved when we vary $p$ , sparsity $k$, or the amount of noise $\sigma$ in the additional blocks of experiments without any loss of accuracy. 
\subsection{Support Vector Machines}
\label{subsec:svm}
In this section, we apply Algorithm  \ref{alg:cutting_plane_sto}
to the problem of support vector machines for classification. In this particular case, the variables are all continuous, and therefore we focus on utilizing the cutting-plane algorithm to solve the NLP subproblem as defined in Algorithm \ref{alg:cutting_plane_sto} since the set of integer variables $\ma{Z}=\emptyset$.  

Specifically, we consider a dataset where we have $(\bm{x}_i,y_i)$, $i\in [N]$  and $y_i \in \{-1,+1\}$, and we wish to create a classifier of the form $f(\bm{x})=\sign(\bm{\theta}^T\bm{x})$. Then the SVM problem with $\ell_2$ regularization can be written as:
\begin{align}
    \min_{\bm{\theta},\epsilon\geq 0} \frac{1}{2}\|\bm{\theta}\|^2+CR(\bm{\theta}) \label{eq:svm},
\end{align}
where $C$ is a user-defined constant controlling the degree of regularization, and:
\begin{equation}
R(\bm{\theta})=\max_{\bm{c} \in \{0,1\}^N}\left\{\frac{1}{N}\sum_{i=1}^N c_i - \frac{1}{N}\sum_{i=1}^N c_iy_i(\bm{w}^T\bm{x}_i)\right\}.
\end{equation}
In this setting  all variables are continuous. We follow the formulation in classical papers of utilizing cutting planes in SVMs (e.g. \cite{franc2008optimized,joachims2006training}) by using cutting planes (i.e. linearizing the objective) only on the $R(\bm{\theta})$ portion of the objective. It is easy to see that $R(\bm{\theta})$ satisfies the data concentration assumption by recognizing the form of $R(\bm{\theta})$ as a maximum over a finite set of expressions in the form of those introduced in Section \ref{subsec:riskmin}, which satisfies Assumption \ref{ass:random} by Lemma \ref{lem:riskmin}. 

For any $\bm{\theta}$, the following expression is a valid subgradient for $R(\bm{\theta})$:
\begin{equation}
    \partial R(\bm{\theta})= \frac{1}{N}\sum_{i=1}^N c_iy_i\bm{x}_i,
\end{equation}
with
\begin{equation*}
    c_i=\begin{cases}
    1, & y_i(\bm{\theta}^T\bm{x}_i)<1, \\
    0, & \text{otherwise.}
    \end{cases}.
\end{equation*}
To conduct the experiment, we  utilize the  forest covertype dataset from the UCI Machine Learning Repository \cite{Dua:2019} that has 580,012 samples with 54 features, and the objective is to predict the type of forest cover. We binarize the objective variable to predict if the forest cover belongs to Class 2 (the most frequent class at $48.9\%$) or not, following a similar treatment in \cite{joachims2006training}. 

We fix $C=10^6$ as defined in Equation (\ref{eq:svm}), and set the tolerance in the cutting plane to $\epsilon=10^{-4}$. Then, we randomly sample two disjoint datasets, each of size $N$, for training and testing from the master dataset. We experiment with different sizes of $N$ to observe the speedup and performance of stochastic cutting planes (SCP) over standard cutting planes. We again set $n=\min\{N, 10\sqrt{N}\}$ as the sampling size for stochastic cutting planes. The results, recording time $T$ and the out-of-sample accuracy on the testing set $ACC$ over 10 random runs are contained in Table \ref{tab:svm}. We further compare to LIBSVM \citep{chang2011libsvm}, an established reference SVM implementation using Sequential Minimal Optimization (SMO), and stochastic gradient descent (SGD), which is widely regarded as the fastest algorithm for solving SVM problems (at the expense of some accuracy). The testing environment has a six-core i7-5820k CPU with 24GB of RAM. Except for SGD, all algorithms are implemented in pure Julia with Julia v1.4.2 to minimize the difference in conditions. For SGD, we utilize the optimized formulation contained in scikit-learn V0.22.2 and Python V3.7.4. 
\begin{table}[H]
	\centering
	\resizebox{0.95\columnwidth}{!}{%
		\begin{tabular}{|l|l||l|c|l|c|l|c|l|c|}
			\hline  \multirow{2}{*}{$\bm{N}$}&\multirow{2}{*}{$\bm{p}$}& \multicolumn{2}{c|}{\textbf{Cutting Planes}} & \multicolumn{2}{c|}{\textbf{SCP}}& \multicolumn{2}{c|}{\textbf{LIBSVM}} & \multicolumn{2}{c|}{\textbf{SGD}}\\\cline{3-10}
			  &    & $\bm{T}$ & \textbf{ACC} & $\bm{T}$ & \textbf{ACC} & $\bm{T}$ & \textbf{ACC}  & $\bm{T}$ & \textbf{ACC}  \\\hline
		 $10^3$ & 54  & 107.6s & $74.9\%$& 7.5s & $72.3\%$ & 0.75s & $74.2\%$ & 0.01s & $71.8\%$\\
			$10^4$ & 54 & 145.7s & $75.5\%$& 10.7s & $74.8\%$& 30.9s & $75.0\%$ & 0.09s & $73.4\%$ \\
					$10^5$ & 54  & 439.2s & $75.9\%$& 32.1s & $75.3\%$& 960.5s & $75.7\%$& 0.45s & $74.3\%$\\\hline

		\end{tabular}
	}
	\caption{Comparison of cutting planes and stochastic cutting planes on the covertype dataset for support vector machines.  }
	\label{tab:svm}
\end{table}
We see that stochastic cutting planes achieve a significant order-of-magnitude speedup compared to the standard cutting plane while maintaining a comparable accuracy with minimum degradation. At the largest scales it is over $13\times $ faster than the standard cutting planes formulation and $30\times$ faster than LIBSVM with only a degradation of $\sim0.5\%$ in out-of-sample accuracy. This reflects the scalability of the stochastic cutting plane method versus other methods in large-scale data. Stochastic gradient descent is still significantly faster, but at the further expense of accuracy of about $1\%$. We note that however, the running time of SCP appears to grow slower (roughly proportional to $\sqrt{n}$) compared to SGD (roughly proportional to $n$). 

We note that at the largest scale with $N=10^5$, each stochastic cutting plane uses $n=10^{3.5}$ samples, and yet it manages to be $3\times $ faster than the full cutting plane algorithm when ran on $10^3$ samples while its accuracy is closer to the cutting plane algorithm when ran on $10^4$ samples. This shows that the randomness of the data samples selected by each cutting plane is key, and the performance of the stochastic cutting plane using a random sample of size $n$ cannot be simply replicated by choosing $n$ samples and running the standard cutting plane algorithm. 
\subsection{The Stochastic Knapsack Problem}
\label{subsec:sskp}
We further apply stochastic cutting planes to SAA versions of stochastic optimization problems. In particular, we test stochastic cutting planes on the classic \emph{static stochastic knapsack problem} (SSKP). A SSKP is a resource allocation problem where a decision maker has to choose a subset of $k$ known alternative projects to take on. For this purpose, a known quantity $q$ of relatively low-cost resource is available, while any further resource required can be obtained at a known cost of $c$ per unit of resource. The amount $W_i$ of resource required by project $i$ is not known exactly but the decision maker has an estimate of the probability distribution. Each project $i$ has an expected net reward of $r_i$. Thus, the optimization problem can be formulated as:
\[\max_{\bm{z} \in \{0,1\}^k} \sum_{i=1}^k r_iz_i-c\E\left[\sum_{i=1}^k W_iz_i-q\right]^+\]
Where $[\cdot]^+:=\max(\cdot,0)$ and the expectation is taken over the probability distribution of $\bm{W}=[W_1,\ldots,W_k]$. This objective function appears in many real-life problems such as airline crew scheduling, shortest path problems, and others, often with further constraints on $\bm{z}$ to limit the optimization over $\bm{z} \in \mathcal{Z} \subset \{0,1\}^k$.

We here consider the SAA version of SSKP where $N$ realizations of $\bm{W}$ are available (either by simulation or historical data) in the form of $[\bm{W}^1,\ldots,\bm{W}^N]$. Then the SAA version of SSKP can be written as:
\begin{equation}
    \max_{\bm{z} \in \{0,1\}^k} \sum_{i=1}^k r_iz_i-\frac{c}{N}\sum_{j=1}^N \left[\sum_{i=1}^k W_i^jz_i-q\right]^+\label{eq:sskp}
\end{equation}
This problem is an integer concave optimization problem, and we can utilize cutting planes on the cost part of the objective:
\begin{equation*}
    C(\bm{z})= \frac{c}{N}\sum_{j=1}^N \left[\sum_{i=1}^k W_i^jz_i-q\right]^+ .
\end{equation*}
$C(\bm{z})$ satisfies the data concentration assumption since it satisfies the form in Lemma \ref{lem:riskmin}. Its derivative is easily calculated as:
\begin{equation*}
    (\nabla C(\bm{z}))_i= \frac{c}{N}\sum_{j=1}^N W_i^j\one\left\{\sum_{i=1}^k W_i^jz_i-q\geq 0\right\}.
\end{equation*}
To generate the specific instances for this problem, we follow the setup in \cite{kleywegt2002sample}. The rewards follow $r_i \sim \Unif[10,20]$ , while $W_i \sim N(\mu_i,\sigma_i)$ where $\mu_i \sim \Unif[20,30]$ and $\sigma_i \in \Unif[5,15]$. For all instances, we set the per unit penalty $c=4$, and set the initial available inventory $q=\max(k,20)$ (the lower bound is so that the optimal solution is not the trivial 0 vector). We investigate different combinations of $(N,k)$ to uncover how the performance of stochastic cutting plane scales with problem size and samples generated. 

We compare stochastic cutting planes (with sampling size $n=\min\{N, 10\sqrt{N}\}$) to the standard cutting plane algorithm, and further with the standard linear reformulation of \eqref{eq:sskp} using auxiliary variables $\bm{x} \in \md{R}^n$ as outlined in \cite{kleywegt2002sample}:
\begin{align}
    &\max_{\bm{z} \in \{0,1\}^k,\bm{x}\geq 0 } \sum_{i=1}^k r_iz_i - \frac{c}{N} \sum_{j=1}^N z_j\label{eq:reformobj}\\
    &s.t. \quad x_j \geq \sum_{i=1}^k W_i^jz_i-q, \;\; \forall j \in [N].\label{eq:reformcons}
\end{align}
The testing environment has a six-core i7-5820k CPU with 24GB of RAM. All algorithms are implemented in pure Julia with Julia v1.4.2 to minimize the difference in conditions, and we utilize Gurobi V9.0.1 as the integer programming solver. We record two statistics: the time taken to solve the formulation ($T$), and the objective value ($Obj$). For ease of comparison, we normalized the objective value of the optimal solution to $100\%$. Each statistic is shown as an average of 20 randomly generated instances:
\begin{table}[H]
	\centering
	\resizebox{0.95\columnwidth}{!}{%
		\begin{tabular}{|l|l||l|c|l|c|l|c|}
			\hline  \multirow{2}{*}{$\bm{N}$}&\multirow{2}{*}{$\bm{k}$}& \multicolumn{2}{c|}{\textbf{Cutting Planes}} & \multicolumn{2}{c|}{\textbf{SCP}}& \multicolumn{2}{c|}{\textbf{Linear Reformulation}} \\\cline{3-8}
			  &    & $\bm{T}$ & \textbf{Obj} & $\bm{T}$ & \textbf{Obj} & $\bm{T}$ & \textbf{Obj}  \\\hline
		 $10^3$ & 10  & 0.33s & $100\%$ & 0.24s & $99.7\%$ & 7.0s & $100\%$ \\
			$10^4$ & 10 & 0.89s & $100\%$ & 0.27s & $99.9\%$ & 70.3s & $100\%$  \\
					$10^5$ & 10  & 3.6s & $100\%$ & 0.35s & $100\%$ & 3760s & $100\%$\\
					$10^6$ & 10  & 61.3s & $100\%$ & 0.97s & $100\%$ & N/A & N/A \\\hline
		 $10^3$ & 20  & 0.42s & $100\%$ & 0.25s & $99.9\%$ & 8.5s & $100\%$\\
			$10^4$ & 20 & 1.1s & $100\%$ & 0.39s & $99.9\%$& 437s & $100\%$ \\
					$10^5$ & 20  & 9.3s & $100\%$  & 0.46s & $100\%$& 7328s & $100\%$\\
					$10^6$ & 20  & 124s & $100\%$ & 1.7s & $100\%$ &N/A  & N/A \\\hline
						 $10^3$ & 50  & 1.3s & $100\%$ & 1.0s & $100\%$ & 6.1s & $100\%$ \\
			$10^4$ & 50 & 5.3s & $100\%$ & 1.1s & $100\%$ & 260s & $100\%$  \\
					$10^5$ & 50  & 51.8s & $100\%$ & 2.5s & $100\%$ & 4005s & $100\%$\\
					$10^6$ & 50  & 706s & $100\%$ & 6.9s & $100\%$ & N/A  & N/A \\\hline

		\end{tabular}
	}
	\caption{Comparison of cutting planes and stochastic cutting planes for SSKP. N/A means that no instance completed running under 10000s. }
	\label{tab:sskp}
\end{table}
We first observe that both the standard and the stochastic cutting plane formulation outperform the linearized reformulation significantly in time needed. This is not surprising as the linear reformulation requires $N$ auxiliary variables to model the $[\cdot]^+$ function, which severely impacts solve time. 

We further notice that stochastic cutting planes is able to achieve a near optimal solution across all combinations of $N$ and $k$ tested. In the largest instances (where $N=10^6$), it was able to almost always return the same solution as the general cutting plane algorithm while spending two orders of magnitude less time. This again demonstrates the significant speedup that stochastic cutting plane provides.

\subsection{Sample Size Selection}
\label{subsec:samplesize}
The key parameter in the method of stochastic cutting planes is the sample size $n$, as it drives both the accuracy in Theorem \ref{theo:convergence}, and the speedup compared to the full cutting plane algorithm. Therefore, in this section, we  explore how the sample size $n$ affects the quality of solutions in various problems.  

Specifically, we  test different combinations of $(N,n)$ on sparse regression with $p=100$, $k=10$ using the same setup as in Section \ref{subsec:sparsereg} over 10 runs each, and record the percentage of cases that the stochastic cutting plane algorithm returned the same solution as the optimal solution. 
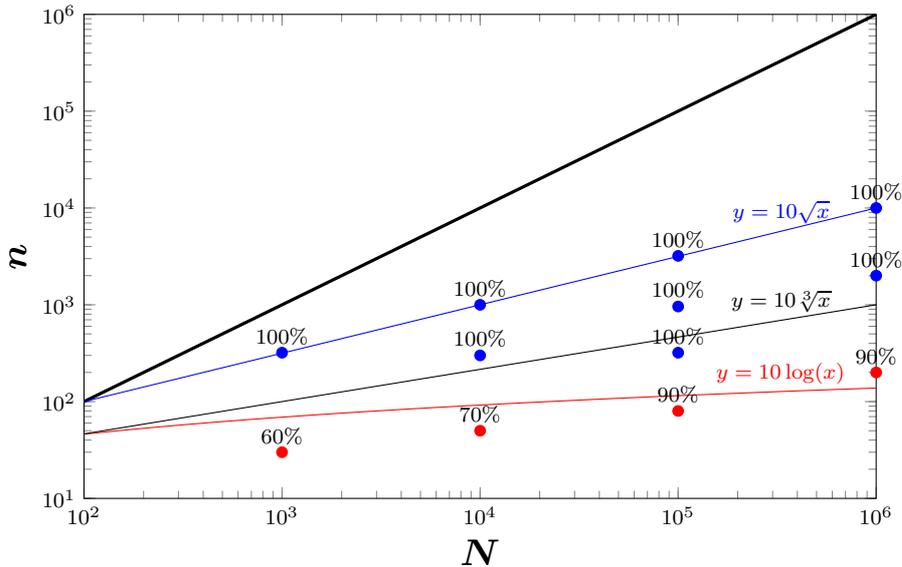
\begin{figure}[htbp]
\centering
\begin{tikzpicture}
\begin{loglogaxis}[width=12cm, height=8cm,
xmin=100,
xmax =1000000,
ymin=10,
ymax =1000000,
xlabel = { $\bm{N}$},
ylabel = {$\bm{n}$},
label style={font=\fontsize{14}{17}}]
    \addplot[
        scatter/classes={a={blue}, b={red}},
        scatter, mark=*, only marks, 
        scatter src=explicit symbolic,
        nodes near coords*={\Label},
        visualization depends on={value \thisrow{label} \as \Label} 
] table [meta=class, row sep = crcr] {
        x y class label\\
        10000 1000 a {$100\%$}\\};
    \addplot[
        scatter/classes={a={blue}, b={red}},
        scatter, mark=*, only marks, 
        scatter src=explicit symbolic,
        nodes near coords*={\Label},
        visualization depends on={value \thisrow{label} \as \Label} 
] table [meta=class, row sep = crcr] {
        x y class label\\
        10000 300 a {$100\%$}\\};
    \addplot[
        scatter/classes={a={blue}, b={red}},
        scatter, mark=*, only marks, 
        scatter src=explicit symbolic,
        nodes near coords*={\Label},
        visualization depends on={value \thisrow{label} \as \Label} 
] table [meta=class, row sep = crcr] {
        x y class label\\
         1000000 10000 a {$100\%$}\\};
    \addplot[
        scatter/classes={a={blue}, b={red}},
        scatter, mark=*, only marks, 
        scatter src=explicit symbolic,
        nodes near coords*={\Label},
        visualization depends on={value \thisrow{label} \as \Label} 
] table [meta=class, row sep = crcr] {
        x y class label\\
         100000 3200 a {$100\%$}\\};
    \addplot[
        scatter/classes={a={blue}, b={red}},
        scatter, mark=*, only marks, 
        scatter src=explicit symbolic,
        nodes near coords*={\Label},
        visualization depends on={value \thisrow{label} \as \Label} 
] table [meta=class, row sep = crcr] {
        x y class label\\
         100000 960 a {$100\%$}\\};
    \addplot[
        scatter/classes={a={blue}, b={red}},
        scatter, mark=*, only marks, 
        scatter src=explicit symbolic,
        nodes near coords*={\Label},
        visualization depends on={value \thisrow{label} \as \Label} 
] table [meta=class, row sep = crcr] {
        x y class label\\
         100000 320 a {$100\%$}\\};
    \addplot[
        scatter/classes={a={blue}, b={red}},
        scatter, mark=*, only marks, 
        scatter src=explicit symbolic,
        nodes near coords*={\Label},
        visualization depends on={value \thisrow{label} \as \Label} 
] table [meta=class, row sep = crcr] {
        x y class label\\
         1000 320 a {$100\%$}\\};
    \addplot[
        scatter/classes={a={blue}, b={red}},
        scatter, mark=*, only marks, 
        scatter src=explicit symbolic,
        nodes near coords*={\Label},
        visualization depends on={value \thisrow{label} \as \Label} 
] table [meta=class, row sep = crcr] {
        x y class label\\
         1000 30 b {$60\%$}\\};
    \addplot[
        scatter/classes={a={blue}, b={red}},
        scatter, mark=*, only marks, 
        scatter src=explicit symbolic,
        nodes near coords*={\Label},
        visualization depends on={value \thisrow{label} \as \Label} 
] table [meta=class, row sep = crcr] {
        x y class label\\
         1000000 2000 a {$100\%$}\\};
    \addplot[
        scatter/classes={a={blue}, b={red}},
        scatter, mark=*, only marks, 
        scatter src=explicit symbolic,
        nodes near coords*={\Label},
        visualization depends on={value \thisrow{label} \as \Label} 
] table [meta=class, row sep = crcr] {
        x y class label\\
         100000 80 b {$90\%$}\\};
\addplot[
        scatter/classes={a={blue}, b={red}},
        scatter, mark=*, only marks, 
        scatter src=explicit symbolic,
        nodes near coords*={\Label},
        visualization depends on={value \thisrow{label} \as \Label} 
] table [meta=class, row sep = crcr] {
        x y class label\\
         1000000 200 b {$90\%$}\\};
\addplot[
        scatter/classes={a={blue}, b={red}},
        scatter, mark=*, only marks, 
        scatter src=explicit symbolic,
        nodes near coords*={\Label},
        visualization depends on={value \thisrow{label} \as \Label} 
] table [meta=class, row sep = crcr] {
        x y class label\\
         10000 50 b {$70\%$}\\};
        \addplot[domain=100:1000000,very thick] {x};   
        \addplot[domain=100:1000000,blue] {10*sqrt(x)} node[above,pos=0.88] {$y=10\sqrt{x}$};
        \addplot[domain=100:1000000,red] {10*ln(x)} node[above,pos=0.88] {$y=10\log(x)$};
        \addplot[domain=100:1000000,black] {10*x^(1/3)} node[above,pos=0.88] {$y=10\sqrt[3]{x}$};

\end{loglogaxis}
\end{tikzpicture}
\caption{Percentage of Runs where the Stochastic Cutting Planes agrees with true solution on different $(N,n)$ combinations for sparse regression with $p=100$ and $k=10$.}
\label{fig:samplesize}
\end{figure}
From Figure \ref{fig:samplesize} we can see that for the task of sparse regression, the sample size $n$ needed so that the stochastic solution always coincides with the true solution seems to grow slower than $O(\sqrt{N})$ but faster than $O(\log(N))$. Therefore, in cases where accuracy of the solution is important, it would seem to be prudent, at least for the task of sparse regression, to select $O(\sqrt{N})$ samples to ensure that the solution does indeed coincide, as we showed in Table \ref{tab:sparsereg}. However, in scenarios where accuracy is less important than scalability, then the practitioner can explore even smaller $n$ such as $O(\sqrt[3]{N})$ as it is clear from the graph that $O(\sqrt{N})$ is a conservative bound as $N \to \infty$. 

We further explore this sample size relationship in the context of the SSKP problem. We follow the same setup as in Section \ref{subsec:sskp}, and set $k=50$. We record the objective value of SCP as a percentage of the optimal solution, averaged over 20 runs. 

\begin{figure}[htbp]
\centering
\begin{tikzpicture}
\begin{loglogaxis}[width=12cm, height=8cm,
xmin=100,
xmax =1000000,
ymin=10,
ymax =1000000,
xlabel = { $\bm{N}$},
ylabel = {$\bm{n}$},
label style={font=\fontsize{14}{17}}]
    \addplot[
        scatter/classes={a={blue}, b={red}},
        scatter, mark=*, only marks, 
        scatter src=explicit symbolic,
        nodes near coords*={\Label},
        visualization depends on={value \thisrow{label} \as \Label} 
] table [meta=class, row sep = crcr] {
        x y class label\\
        1000 300 a {$100\%$}\\};
    \addplot[
        scatter/classes={a={blue}, b={red}},
        scatter, mark=*, only marks, 
        scatter src=explicit symbolic,
        nodes near coords*={\Label},
        visualization depends on={value \thisrow{label} \as \Label} 
] table [meta=class, row sep = crcr] {
        x y class label\\
        1000 100 a {$99.9\%$}\\};
    \addplot[
        scatter/classes={a={blue}, b={red}},
        scatter, mark=*, only marks, 
        scatter src=explicit symbolic,
        nodes near coords*={\Label},
        visualization depends on={value \thisrow{label} \as \Label} 
] table [meta=class, row sep = crcr] {
        x y class label\\
        1000 30 b {$99.5\%$}\\};
    \addplot[
        scatter/classes={a={blue}, b={red}},
        scatter, mark=*, only marks, 
        scatter src=explicit symbolic,
        nodes near coords*={\Label},
        visualization depends on={value \thisrow{label} \as \Label} 
] table [meta=class, row sep = crcr] {
        x y class label\\
        10000 1000 a {$99.99\%$}\\};
    \addplot[
        scatter/classes={a={blue}, b={red}},
        scatter, mark=*, only marks, 
        scatter src=explicit symbolic,
        nodes near coords*={\Label},
        visualization depends on={value \thisrow{label} \as \Label} 
] table [meta=class, row sep = crcr] {
        x y class label\\
        10000 100 b {$99.78\%$}\\};
    \addplot[
        scatter/classes={a={blue}, b={red}},
        scatter, mark=*, only marks, 
        scatter src=explicit symbolic,
        nodes near coords*={\Label},
        visualization depends on={value \thisrow{label} \as \Label} 
] table [meta=class, row sep = crcr] {
        x y class label\\
        10000 10 b {$98.9\%$}\\};
    \addplot[
        scatter/classes={a={blue}, b={red}},
        scatter, mark=*, only marks, 
        scatter src=explicit symbolic,
        nodes near coords*={\Label},
        visualization depends on={value \thisrow{label} \as \Label} 
] table [meta=class, row sep = crcr] {
        x y class label\\
         100000 3200 a {$99.99\%$}\\};
    \addplot[
        scatter/classes={a={blue}, b={red}},
        scatter, mark=*, only marks, 
        scatter src=explicit symbolic,
        nodes near coords*={\Label},
        visualization depends on={value \thisrow{label} \as \Label} 
] table [meta=class, row sep = crcr] {
        x y class label\\
         100000 320 a {$99.91\%$}\\};
    \addplot[
        scatter/classes={a={blue}, b={red}},
        scatter, mark=*, only marks, 
        scatter src=explicit symbolic,
        nodes near coords*={\Label},
        visualization depends on={value \thisrow{label} \as \Label} 
] table [meta=class, row sep = crcr] {
        x y class label\\
         100000 32 b {$99.62\%$}\\};        
    \addplot[
        scatter/classes={a={blue}, b={red}},
        scatter, mark=*, only marks, 
        scatter src=explicit symbolic,
        nodes near coords*={\Label},
        visualization depends on={value \thisrow{label} \as \Label} 
] table [meta=class, row sep = crcr] {
        x y class label\\
         1000000 100 a {$99.90\%$}\\};    
    \addplot[
        scatter/classes={a={blue}, b={red}},
        scatter, mark=*, only marks, 
        scatter src=explicit symbolic,
        nodes near coords*={\Label},
        visualization depends on={value \thisrow{label} \as \Label} 
] table [meta=class, row sep = crcr] {
        x y class label\\
         1000000 1000 a {$99.99\%$}\\};  
    \addplot[
        scatter/classes={a={blue}, b={red}},
        scatter, mark=*, only marks, 
        scatter src=explicit symbolic,
        nodes near coords*={\Label},
        visualization depends on={value \thisrow{label} \as \Label} 
] table [meta=class, row sep = crcr] {
        x y class label\\
         1000000 10000 a {$99.99\%$}\\};  
        \addplot[domain=100:1000000,very thick] {x};   
        \addplot[domain=100:1000000,blue] {10*sqrt(x)} node[above,pos=0.88] {$y=10\sqrt{x}$};
        \addplot[domain=100:1000000,red] {10*ln(x)} node[above,pos=0.88] {$y=10\log(x)$};
        \addplot[domain=100:1000000,black] {10*x^(1/3)} node[above,pos=0.88] {$y=10\sqrt[3]{x}$};

\end{loglogaxis}
\end{tikzpicture}
\caption{Average performance of the solution reached by Stochastic Cutting Planes as a percentage of true optimal objective on different $(N,n)$ combinations for SSKP with $k=50$. Nodes with $<99.9\%$ of optimality are colored red.}
\label{fig:samplesize_sskp}
\end{figure}
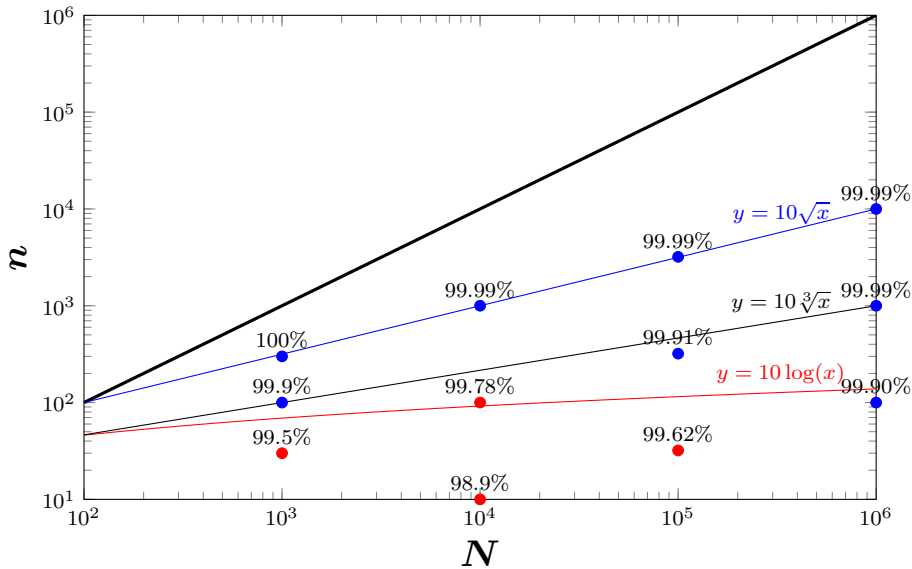
Similar to sparse regression, we see from Figure \ref{fig:samplesize_sskp} that the chosen sample size $n=10\sqrt{N}$ is a quite conservative and the average objective across 20 runs is always within $0.01\%$ of the optimal solution. We further see that even at scales of $n=10\log(N)$, the objective value of the solution reached by SCP is always within $0.5\%$ of the true solution. Therefore, in problems where absolute optimality is not necessary, further running time gains (compared to what is demonstrated in this work) could be achieved by using a smaller sample size. 

\section{Conclusions}
We have presented stochastic cutting planes for data-driven convex MINLPs. Our theoretical results show that the proposed methodology is able to obtain the optimal solution with exponentially vanishing failure probability. Our numerical experiments demonstrate that on a wide variety of problems, the stochastic cutting planes algorithm is able to achieve an order-of-magnitude speedup compared to the standard cutting planes algorithm without a significant sacrifice in accuracy, and is competitive with state-of-the-art algorithms. 

\bibliographystyle{spbasic} 
\bibliography{SCPbib} 

\clearpage
\appendix

\setcounter{table}{0}
\renewcommand{\thetable}{A\arabic{table}}
\setcounter{figure}{0}
\renewcommand{\thefigure}{A\arabic{figure}}
\setcounter{equation}{0}
\renewcommand{\theequation}{A\arabic{equation}}
\section{Proof of Lemma \ref{lem:sparsereg}}
\label{app:lem_sparsereg_proof}
We only prove the lemma  for the function $f(\bm{z},\bm{\theta}; [N])$ as the derivative $\nabla f$ follows the exact same procedure just with additional notation. 

First, by the matrix inversion lemma \citep{matrixinv}, for any feasible solution $\bm{z} \in \ma{Z}$, we can rewrite $f(\bm{z},\bm{\theta}; [N])$ as:
\begin{align}
    f(\bm{z},\bm{\theta}; [N])&=\frac{1}{N}\bm{y}^T\left(I_N-\bm{X}_{\bm{z}}\left(\frac{I_k}{\gamma} + \bm{X}_{\bm{z}}^T\bm{X}_{\bm{z}}\right)^{-1}\bm{X}_{\bm{z}}\right)\bm{y}\nonumber \\&=\frac{1}{N}\bm{y}^T\bm{y}-\frac{1}{N}\bm{y}^T\bm{X}_{\bm{z}}\left(\frac{I_k}{\gamma} + \bm{X}_{\bm{z}}^T\bm{X}_{\bm{z}}\right)^{-1}\bm{X}_{\bm{z}}\bm{y}\nonumber \\&=\frac{1}{N}\bm{y}^T\bm{y}-\frac{1}{N}\bm{y}^T\bm{V}\left(\frac{I_k}{\gamma} + \bm{V}^T\bm{V}\right)^{-1}\bm{V}\bm{y}, \label{eq:fullsparsereg}
\end{align}
where $\bm{V}=\bm{X}_{\bm{z}} \in \md{R}^{N \times k}$ is the subset of the design matrix $\bm{X}$ with only the columns $\{i\mid z_i=1\}$. For simplicity of notation, we use $\bm{V}$ to suppress the dependence on $\bm{z}$.

We let $S$ be a subset of $[N]$ of size $n$ (we suppress the dependence of $S$ on $n$ to reduce notation complexity). Then, we  have, similar to Equation (\ref{eq:fullsparsereg}):
\begin{equation}
     f(\bm{z},\bm{\theta}; S)=\frac{1}{n}\bm{y}_S^T\bm{y}_S-\frac{1}{n}\bm{y}_S^T\bm{V}_S\left(\frac{I_k}{\gamma} + \bm{V}_S^T\bm{V}_S\right)^{-1}\bm{V}_S\bm{y}_S,
\end{equation}
where $\bm{y}_S$ is the subvector of $\bm{y}$ with only the entries in set $S$, and similar for $\bm{V}$. Therefore, what we need to prove is that:
\begin{align}
    &\P\left(\left|\frac{1}{N}\bm{y}^T\bm{y}-\frac{1}{N}\bm{y}^T\bm{V}\left(\frac{I_k}{\gamma} + \bm{V}^T\bm{V}\right)^{-1}\bm{V}\bm{y}-\frac{1}{n}\bm{y}_S^T\bm{y}_S+\frac{1}{n}\bm{y}_S^T\bm{V}_S\left(\frac{I_k}{\gamma} + \bm{V}_S^T\bm{V}_S\right)^{-1}\bm{V}_S\bm{y}_S\right|\leq M \sqrt{\frac{\log(\frac{1}{\epsilon})}{n}}\right)\nonumber \\&\geq 1- \epsilon. \label{eq:concentration_main}
\end{align}

We now introduce and prove lemmas that allow us to prove Equation (\ref{eq:concentration_main}):
\begin{lem}
\label{lem:matvecconcentration}
Let $S$ be a random subset of $[N]$ of size $n$. Then, we have
\begin{align}
    \md{P}\left(\left|\frac{\bm{y}_S^T\bm{y}_S}{n}-\frac{\bm{y}^T\bm{y}}{N}\right|\leq \sqrt{\frac{A\log(\frac{1}{\epsilon})}{n}}\right)&\geq 1- \epsilon, \label{eq:matvec_dev_bound1}\\
    \md{P}\left(\left\|\frac{\bm{y}_S^T\bm{V}_S}{n}-\frac{\bm{y}^T\bm{V}}{N}\right\|\leq \sqrt{\frac{B\log(\frac{1}{\epsilon})}{n}}\right)&\geq 1- \epsilon. \label{eq:matvec_dev_bound2}   
\end{align}
\end{lem}
\proof{Proof:}
Notice that we have:
\begin{alignat*}{3}
    \bm{y}^T\bm{y}&=\sum_{i=1}^N y_i^2,\qquad  &\qquad  \bm{y}^T\bm{V}&=\sum_{i=1}^N y_i\bm{v}_i,\\
    \bm{y}_S^T\bm{y}_S&=\sum_{i\in S} y_i^2,\qquad  &\qquad  \bm{y}_S^T\bm{V}_S&=\sum_{i\in S} y_i\bm{v}_i.
\end{alignat*}
Therefore, if we treat $y_1^2,\ldots y_N^2$ as a finite population, then $\bm{y}_S^T\bm{y}_S$ is a random sample of $n$ points drawn without replacement from that set, and similarly for $\bm{y}_S^T\bm{V}_S$. Then, the required inequality follows directly from Hoeffding's Inequality in the case without replacement, stated below:
\begin{prop}[Hoeffding's Inequality]
\label{prop:hoeffding}
Let $\mathcal{X}=(x_1,\ldots,x_n)$ be a finite population of $N$ points and $X_1,\ldots, X_n$ be a random sample drawn without replacement from $\bm{X}$. Let:
\begin{equation*}
    a = \min_{1\leq i \leq n} x_i, \qquad \text{and} \qquad b = \max_{1\leq i \leq n} x_i.
\end{equation*}
Then, for all $\epsilon>0$, we have:
\begin{equation}
    \md{P}\left(\Big|\frac{\sum_{i=1}^n X_i}{n} -\mu \Big|\geq \epsilon\right) \leq 2\exp\left(-\frac{2n\epsilon^2}{(b-a)^2}\right).
\end{equation}
\end{prop}
For proof of this proposition, see for example \cite{boucheron2013concentration}. \QED
\endproof
Then, we have the following lemma:
\begin{lem}
\label{lem:matrix_ineq}
\begin{equation}
    \P \left(\left\|\left(\frac{I_k}{n\gamma} + \frac{\bm{V}^T_S\bm{V}_S}{n}\right)^{-1}-\left(\frac{I_k}{N\gamma} + \frac{\bm{V}^T\bm{V}}{N}\right)^{-1}\right\|\leq \sqrt{\frac{H\log(\frac{1}{\epsilon})}{n}}\right)\geq 1- \epsilon . 
\end{equation}
\end{lem}
To prove this, we first introduce a matrix analog of the well-known Chernoff bound, the proof of which can be found in \cite{tropp2012user}:
\begin{lem}
    \label{lem:matrix_chernoff}
    Let $\mathcal{X} \in \md{R}^{k \times k}$ be a finite set of positive-semidefinite matrices, and suppose that:
    \[\max_{\bm{X} \in \mathcal{X}} \lambda_{\max}(\bm{X})\leq D,\]
    where $\lambda_{\min}/\lambda_{\max}$is the minimum/maximum eigenvalue function. Sample $\{\bm{X}_1,\ldots, \bm{X}_\ell\}$ uniformly at random without replacement. Let
    \[\mu_{\min}:=\ell \cdot \lambda_{\min}(\md{E} \bm{X}_1),\qquad \mu_{\max}:=\ell \cdot \lambda_{\max}(\md{E} \bm{X}_1).\]
    Then
    \begin{align*}
        \md{P}\left\{\lambda_{\min} \left(\sum_j \bm{X}_j\right)\leq (1-\delta) \mu_{\min} \right\}&\leq k \cdot \exp\left(\frac{-\delta^2\mu_{\min}}{4D}\right), \quad \text{for } \delta \in [0,1),\\
        \md{P}\left\{\lambda_{\max} \left(\sum_j \bm{X}_j\right)\leq (1+\delta) \mu_{\max} \right\}&\leq k \cdot \exp\left(\frac{-\delta^2\mu_{\max}}{4D}\right), \quad \text{for } \delta \geq 0. 
    \end{align*}
\end{lem}
Now we proceed with the proof.
\proof{Proof of Lemma \ref{lem:matrix_ineq}:}
First, we  do a QR decomposition of the matrix $\bm{V}=\bm{Q}\bm{R}$ where $\bm{Q} \in \md{R}^{N \times k}$ has orthogonal columns and $\bm{R} \in \md{R}^{k \times k}$ is upper triangular. We further normalize the decomposition so we have that $\frac{\bm{Q}^T\bm{Q}}{N}=\bm{I}_k$. Then, we have  that $\bm{V}_S=\bm{Q}_S\bm{R}$, and we note that:
\begin{align*}
 \bm{Q}^T\bm{Q}&=\sum_{i=1}^N \bm{q}_i\bm{q}_i^T,\\
 \bm{Q}^T_S\bm{Q}_S&=\sum_{i\in S} \bm{q}_i\bm{q}_i^T,
\end{align*}
where $\bm{q}_i\bm{q}_i^T \in \md{R}^{k \times k}$ rank-one positive semi-definite matrices. Therefore, we can take $\bm{Q}^T_S\bm{Q}_S$ as a random sample of size $n$ from the set $\mathcal{X}=\{\bm{q}_i\bm{q}_i^T\}_{i\in [N]}$, which satisfies the conditions in Lemma \ref{lem:matrix_chernoff} with $D=O(k)$. Therefore, we use Lemma \ref{lem:matrix_chernoff} to bound $\bm{Q}^T_S\bm{Q}_S$ that:
\begin{align*}
    \md{P}\left\{\lambda_{\min} \left(\frac{\bm{Q}^T_S\bm{Q}_S}{n}\right)\leq (1-\delta)\right\}&\leq k \cdot \exp\left(\frac{-\delta^2n}{kD}\right),\\
        \md{P}\left\{\lambda_{\max} \left(\frac{\bm{Q}^T_S\bm{Q}_S}{n}\right)\geq (1+\delta)\right\}&\leq k \cdot \exp\left(\frac{-\delta^2n}{kD}\right),
\end{align*}
for some absolute constant $D$. Some rearrangement gives:
\begin{equation}
\label{eq:matrix_eig_ineq}
    \md{P}\left\{\lambda_{\min} \left(\frac{\bm{Q}_S^T\bm{Q}_S}{n}\right)\geq 1-\sqrt{\frac{kD\log\left(\frac{2k}{\epsilon}\right)}{n}} \;\;  \text{and} \;\; \lambda_{\max} \left(\frac{\bm{Q}_S^T\bm{Q}_S}{n}\right)\leq 1+\sqrt{\frac{kD\log\left(\frac{2k}{\epsilon}\right)}{n}}\right\}\geq 1 - \epsilon.
\end{equation}
Now since $\frac{\bm{Q}^T\bm{Q}}{N}=\bm{I}_k$, we have:
\begin{equation}
\label{eq:matrix_eig_ineq2}
    \lambda_{\min}\left(\frac{\bm{Q}^T\bm{Q}}{N}\right)=\lambda_{\max}\left(\frac{\bm{Q}^T\bm{Q}}{N}\right)=1.
\end{equation}
Combining equation (\ref{eq:matrix_eig_ineq2})
 and (\ref{eq:matrix_eig_ineq}) we obtain
\begin{equation}
    \md{P}\left\{\left\|\frac{\bm{Q}_S^T\bm{Q}_S}{n}-\frac{\bm{Q}^T\bm{Q}}{N}\right\|\leq \sqrt{\frac{kD\log\left(\frac{2k}{\epsilon}\right)}{n}}\right\}\geq 1 - \epsilon.
\end{equation}
Thus, we have 
\begin{align*}
    &\md{P}\left\{\left\|\frac{\bm{V}_S^T\bm{V}_S}{n}-\frac{\bm{V}^T\bm{V}}{N}\right\|\leq\|\bm{R}\|^2 \sqrt{\frac{kD\log\left(\frac{2k}{\epsilon}\right)}{n}}\right\}\\&\geq\md{P}\left\{\left\|\frac{\bm{Q}_S^T\bm{Q}_S}{n}-\frac{\bm{Q}^T\bm{Q}}{N}\right\|\leq \sqrt{\frac{kD\log\left(\frac{2k}{\epsilon}\right)}{n}}\right\}\geq 1 - \epsilon.
\end{align*}
We let $D'=\log(2)\|\bm{R}\|^4D$ and obtain
\begin{equation}
    \md{P}\left\{\left\|\frac{\bm{V}_S^T\bm{V}_S}{n}-\frac{\bm{V}^T\bm{V}}{N}\right\|\leq \sqrt{\frac{kD'\log\left(\frac{k}{\epsilon}\right)}{n}}\right\}\geq 1 - \epsilon.\label{eq:matrix_eig_ineq3}
\end{equation}
From this, we can easily see that there exists a constant $G$ such that:
\begin{equation}
    \md{P}\left\{\left\|\left(\frac{\bm{I}_k}{n\gamma} + \frac{\bm{V}_S^T\bm{V}_S}{n}\right)-\left(\frac{\bm{I}_k}{N\gamma}+\frac{\bm{V}^T\bm{V}}{N}\right)\right\|\leq \sqrt{\frac{kG\log\left(\frac{k}{\epsilon}\right)}{n}}\right\}\geq 1 - \epsilon,\label{eq:matrix_eig_ineq4}
\end{equation}
as $\frac{\bm{I}_k}{n\gamma}$ and $\frac{\bm{I}_k}{N\gamma}$ are of order of magnitude $O(\frac{1}{n})$. Now let us introduce another lemma from matrix perturbation theory (for proof, see e.g. \cite{stewart1990matrix}):
\begin{lem}
    \label{lem:matrix_inv}
    Let $\bm{A},\bm{B}$ be invertible matrices. Then, we have the following bound:
    \begin{equation}
        \|\bm{A}^{-1}-\bm{B}^{-1}\|\leq \|\bm{A}^{-1}\|\|\bm{B}^{-1}\|\|\bm{A}-\bm{B}\|.
    \end{equation}
\end{lem}
Thus, now let $H=G\left\|\left(\frac{\bm{I}_k}{n\gamma} + \frac{\bm{V}_S^T\bm{V}_S}{n}\right)^{-1}\right\|^2\left\|\left(\frac{\bm{I}_k}{N\gamma} + \frac{\bm{V}^T\bm{V}}{N}\right)^{-1}\right\|^2$, be a constant. Then we have, using Lemma \ref{lem:matrix_inv} and Equation (\ref{eq:matrix_eig_ineq4}):
\begin{align*}
     &\md{P}\left\{\left\|\left(\frac{\bm{I}_k}{n\gamma} + \frac{\bm{V}_S^T\bm{V}_S}{n}\right)^{-1}-\left(\frac{\bm{I}_k}{N\gamma}+\frac{\bm{V}^T\bm{V}}{N}\right)^{-1}\right\|\leq \sqrt{\frac{kH\log\left(\frac{k}{\epsilon}\right)}{n}}\right\}\\&\geq \md{P}\left\{\left\|\left(\frac{\bm{I}_k}{n\gamma} + \frac{\bm{V}_S^T\bm{V}_S}{n}\right)-\left(\frac{\bm{I}_k}{N\gamma}+\frac{\bm{V}^T\bm{V}}{N}\right)\right\|\left\|\left(\frac{\bm{I}_k}{n\gamma} + \frac{\bm{V}_S^T\bm{V}_S}{n}\right)^{-1}\right\|\left\|\left(\frac{\bm{I}_k}{N\gamma} + \frac{\bm{V}^T\bm{V}}{N}\right)^{-1}\right\|\leq \sqrt{\frac{kH\log\left(\frac{k}{\epsilon}\right)}{n}}\right\}\\&=\md{P}\left\{\left\|\left(\frac{\bm{I}_k}{n\gamma} + \frac{\bm{V}_S^T\bm{V}_S}{n}\right)-\left(\frac{\bm{I}_k}{N\gamma}+\frac{\bm{V}^T\bm{V}}{N}\right)\right\|\leq \sqrt{\frac{kG\log\left(\frac{k}{\epsilon}\right)}{n}}\right\}\\& \geq 1-\epsilon.
\end{align*}
As required. We then absorb $H'=k\log(k)H$ as $k$ is a dimension independent from the variable dimension $p$ and the sample dimension $n$. 
\QED
\endproof
With Lemma  \ref{lem:matvecconcentration} and \ref{lem:matrix_ineq}, we are now ready to prove the main result in Equation (\ref{eq:concentration_main}). First, let $M$ be a constant such that:
\[M\geq \max\left\{\sqrt{4A}, 36\left\|\left(\frac{I_k}{N\gamma} + \frac{\bm{V}^T\bm{V}}{N}\right)^{-1}\frac{\bm{V}\bm{y}}{N}\right\|^2B, 36\left\|\frac{\bm{y}_S^T\bm{V}_S}{n}\right\|\left\|\frac{\bm{V}\bm{y}}{N}\right\|^2H,36\left\|\left(\frac{I_k}{n\gamma} + \frac{\bm{V}_S^T\bm{V}_S}{n}\right)^{-1}\frac{\bm{V}_S\bm{y}_S}{n}\right\|^2B\right\}.\]
We introduce the Frechet inequalities, which state that:
    $$\P(A_1\cap A_2 \cap \ldots  \cap A_n) \geq  \P(A_1) + \P(A_2) + \ldots + \P(A_n) -(n-1).$$
Then we have, by Frechet's inequalities:
\begin{align*}
&\P\left(\left|\frac{1}{N}\bm{y}^T\bm{y}-\frac{1}{N}\bm{y}^T\bm{V}\left(\frac{I_k}{\gamma} + \bm{V}^T\bm{V}\right)^{-1}\bm{V}\bm{y}-\frac{1}{n}\bm{y}_S^T\bm{y}_S+\frac{1}{n}\bm{y}_S^T\bm{V}_S\left(\frac{I_k}{\gamma} + \bm{V}_S^T\bm{V}_S\right)^{-1}\bm{V}_S\bm{y}_S\right|\leq M \sqrt{\frac{\log(\frac{1}{\epsilon})}{n}}\right)\\&\geq \P\left(\left|\frac{1}{N}\bm{y}^T\bm{y}-\frac{1}{n}\bm{y}_S^T\bm{y}_S\right|\leq \frac{M}{2} \sqrt{\frac{\log(\frac{1}{\epsilon})}{n}}\right) + \\&+\P\left(\left|\frac{1}{N}\bm{y}^T\bm{V}\left(\frac{I_k}{\gamma} + \bm{V}^T\bm{V}\right)^{-1}\bm{V}\bm{y}-\frac{1}{n}\bm{y}_S^T\bm{V}_S\left(\frac{I_k}{\gamma} + \bm{V}_S^T\bm{V}_S\right)^{-1}\bm{V}_S\bm{y}_S\right|\leq \frac{M}{2} \sqrt{\frac{\log(\frac{1}{\epsilon})}{n}}\right)-1.&\intertext{Again, by Frechet inequalities, we can further breakdown the expression above to:}&\geq \P\left(\left|\frac{1}{N}\bm{y}^T\bm{y}-\frac{1}{n}\bm{y}_S^T\bm{y}_S\right|\leq \frac{M}{2} \sqrt{\frac{\log(\frac{1}{\epsilon})}{n}}\right) \\&+\P\left(\left|\left(\frac{1}{N}\bm{y}^T\bm{V}-\frac{1}{n}\bm{y}_S^T\bm{V}_S\right)\left(\frac{I_k}{N\gamma} + \frac{\bm{V}^T\bm{V}}{N}\right)^{-1}\frac{\bm{V}\bm{y}}{N}\right|\leq \frac{M}{6} \sqrt{\frac{\log(\frac{1}{\epsilon})}{n}}\right)\\&+\P\left(\left|\frac{1}{n}\bm{y}_S^T\bm{V}_S\left(\left(\frac{I_k}{N\gamma} + \frac{\bm{V}^T\bm{V}}{N}\right)^{-1}-\left(\frac{I_k}{n\gamma} + \frac{\bm{V}_S^T\bm{V}_S}{n}\right)^{-1}\right)\frac{\bm{V}\bm{y}}{N}\right|\leq \frac{M}{6} \sqrt{\frac{\log(\frac{1}{\epsilon})}{n}}\right)\\&+\P\left(\left|\frac{1}{n}\bm{y}_S^T\bm{V}_S\left(\frac{I_k}{n\gamma} + \frac{\bm{V}_S^T\bm{V}_S}{n}\right)^{-1}\left(\frac{\bm{V}\bm{y}}{n}-\frac{\bm{V}_S\bm{y}_S}{n}\right)\right|\leq \frac{M}{6} \sqrt{\frac{\log(\frac{1}{\epsilon})}{n}}\right)-3&\intertext{By the definition of $M$, we can simplify the expressions above to:}&\geq \P\left(\left|\frac{1}{N}\bm{y}^T\bm{y}-\frac{1}{n}\bm{y}_S^T\bm{y}_S\right|\leq  \sqrt{\frac{A\log(\frac{1}{\epsilon})}{n}}\right) \\&+\P\left(\left\|\frac{1}{N}\bm{y}^T\bm{V}-\frac{1}{n}\bm{y}_S^T\bm{V}_S\right\|\leq  \sqrt{\frac{B\log(\frac{1}{\epsilon})}{n}}\right)\\&+\P\left(\left\|\left(\frac{I_k}{N\gamma} + \frac{\bm{V}^T\bm{V}}{N}\right)^{-1}-\left(\frac{I_k}{n\gamma} + \frac{\bm{V}_S^T\bm{V}_S}{n}\right)^{-1}\right\|\leq  \sqrt{\frac{H\log(\frac{1}{\epsilon})}{n}}\right)\\&+\P\left(\left\|\frac{\bm{V}\bm{y}}{n}-\frac{\bm{V}_S\bm{y}_S}{n}\right\|\leq  \sqrt{\frac{B\log(\frac{1}{\epsilon})}{n}}\right)-3.&\intertext{Now we apply Lemma \ref{lem:matvecconcentration} and \ref{lem:matrix_ineq} to evaluate these expressions:}&\geq 4-4\epsilon-3 \\& = 1-4\epsilon.
\end{align*}
Thus, taking $\epsilon'=4\epsilon$, and $M'=\frac{M}{\log(2)}$ gives the result we want. 
\vfill
\pagebreak
\section{Proof of Theorem \ref{theo:convergence}}
\label{app:theo_converge_proof}

\subsection{Proof of Finiteness}

We first prove that the outer approximation algorithm indeed converges in finite number of steps $C$ with probability 1. To do so, we note that the set of feasible solutions for the integer variables $\ma{Z}$ is finite, and we prove that the series of feasible integer solutions $\bm{z}_1,\bm{z}_2,\cdots$ created during Algorithm \ref{alg:cutting_plane_sto} never repeats indefinitely, with probability 1.

First, we show that if we have $\bm{z}_a=\bm{z}_b$ for some $a+1\leq b$ where we have not terminated at iteration $b$, then we must have $f(\bm{z}_b,\bm{\theta}_b;[S_b])>f(\bm{z}_a,\bm{\theta}_a;S_a])$.  At iteration $a$, the following cutting plane was added:
\begin{equation}
\eta\geq f(\bm{z}_a,\bm{\theta}_a;[S_a])+\nabla f(\bm{z}_a,\bm{\theta}_a;[S_a])^T\binom{\bm{z}-\bm{z}_a}{\bm{\theta} - \bm{\theta}_a}. \label{eq:cutplane_ineq0}
\end{equation}
Therefore, using (\ref{eq:cutplane_ineq0}) at iteration $b-1\geq a$ (note $\bm{z}_b, \eta_b$ is the optimal solution in the $b-1$th iteration, so it must satisfy the inequality in (\ref{eq:cutplane_ineq0})) gives:
\begin{align}
\eta_b&\geq f(\bm{z}_a,\bm{\theta}_a;[S_a])+\nabla f(\bm{z}_a,\bm{\theta}_a;[S_a])^T\binom{\bm{z}_b-\bm{z}_a}{\bm{\theta}_b - \bm{\theta}_a}\\&=f(\bm{z}_a,\bm{\theta}_a;[S_a])+\nabla_{\bm{\theta}} f(\bm{z}_a,\bm{\theta}_a;[S_a])^T(\bm{\theta}_b - \bm{\theta}_a)\\&\geq f(\bm{z}_a,\bm{\theta}_a;[S_a]).\label{eq:cutplane_ineq}
\end{align}
The last inequality follows because $\bm{\theta}_a$ is the optimal solution to $\NLP(\bm{\theta}_a;[S_a])$, so no feasible direction can exist at $\bm{\theta}_a$. Now consider the cutting plane problem at iteration $b$, where the following cutting plane was added:
\begin{equation}
\eta\geq f(\bm{z}_b,\bm{\theta}_b;[S_b])+\nabla f(\bm{z}_b,\bm{\theta}_b;[S_b])^T\binom{\bm{z}-\bm{z}_b}{\bm{\theta} - \bm{\theta}_b}. \label{eq:cutplane_ineq1}
\end{equation}
Assume that we have  $f(\bm{z}_b,\bm{\theta}_b;[S_b])\leq f(\bm{z}_a,\bm{\theta}_a;[S_a])$. Then $\eta=\eta_b$, $\bm{z}=\bm{z}_b$ is a feasible solution for the cutting plane problem at iteration $b$ since it satisfies the newly added constraint: $\eta_b\geq f(\bm{z}_a,\bm{\theta}_a;[S_a])\geq f(\bm{z}_b,\bm{\theta}_b,[S_b])+\nabla f(\bm{z}_b,\bm{\theta}_b;[S_b])^T\binom{\bm{z}_b-\bm{z}_b}{\bm{\theta}_b - \bm{\theta}_b}=f(\bm{z}_b,\bm{\theta}_b;[S_b])$. By construction of the algorithm, we therefore satisfy the termination condition, so the algorithm terminates, a contradiction. Therefore, we must have $f(\bm{z}_b,\bm{\theta}_b;[S_b])>f(\bm{z}_a,\bm{\theta}_a;[S_a])$.

Then consider Algorithm \ref{alg:cutting_plane_sto} after $Q=|\ma{Z}| (q-1)+1$ iterations where $q>0$ is a positive integer, with the corresponding series of feasible integer solutions being $\bm{z}_1,\cdots,\bm{z}_Q$. Then since there are only $|\ma{Z}|$ feasible solutions, by the pigeonhole principle, there must be $\bm{z} \in \ma{Z}$ that has appeared at least $q$ times in the series, and denote these as $\bm{z}_{l_1},\cdots,\bm{z}_{l_q}$ . By the observation above, we know that if Algorithm \ref{alg:cutting_plane_sto} has not terminated, we must have:
\[f(\bm{z}_{l_1},\bm{\theta}_{l_1};[S_{l_1}])<f(\bm{z}_{l_2},\bm{\theta}_{l_2};[S_{l_2}])<\cdots<f(\bm{z}_{l_q},\bm{\theta}_{l_q};[S_{l_q}])\]
Now note that $\bm{z}_{l_1}=\bm{z}_{l_2}=\cdots=\bm{z}_{l_q}$. Given $\bm{z}_t$, $f(\bm{z}_t,\bm{\theta}_t;[S_t])$ is a random variable that only depends on $S_t$ (by construction of the algorithm,  $\bm{\theta}_t$ is the optimal solution to $\NLP(\bm{z}_t,[S_t])$. Therefore, we have that $f(\bm{z}_{l_1},\bm{\theta}_{l_1};[S_{l_1}]),\cdots,f(\bm{z}_{l_q},\bm{\theta}_{l_q};[S_{l_q}])$ are $q$ independent realizations (since $S_{l_1},S_{l_2},\cdots$ are independent) of an identically distributed random variable, and thus by standard probability results (see e.g. \cite{feller2008introduction}) we have that:
\[\P(f(\bm{z}_{l_1},\bm{\theta}_{l_1};[S_{l_1}])<f(\bm{z}_{l_2},\bm{\theta}_{l_2};[S_{l_2}])<\cdots<f(\bm{z}_{l_q},\bm{\theta}_{l_q};[S_{l_q}]))\leq \frac{1}{q!}\]
The result follows by noting that $f(\bm{z}_{l_1},\bm{\theta}_{l_1};[S_{l_1}]),\cdots,f(\bm{z}_{l_q},\bm{\theta}_{l_q};[S_{l_q}])$ are exchangable random variables so any ordering is equally likely, and only at most 1 out of all $q!$ permutations of $f(\bm{z}_{l_1},\bm{\theta}_{l_1};[S_{l_1}]),\cdots,f(\bm{z}_{l_q},\bm{\theta}_{l_q};[S_{l_q}])$ satisfy the monotone relationship. Then we have that:
\begin{align*}
&\P(\text{Algorithm has not terminated after $Q$ iterations}\mid \bm{z}_1,\cdots \bm{z}_Q)&\\&\leq\P(f(\bm{z}_{l_1},\bm{\theta}_{l_1};[S_{l_1}])<f(\bm{z}_{l_2},\bm{\theta}_{l_2};[S_{l_2}])<\cdots<f(\bm{z}_{l_q},\bm{\theta}_{l_q};[S_{l_q}]))\\&\leq \frac{1}{q!}
\end{align*}
Therefore, by taking expectations on both sides, we have:
\[\P(\text{Algorithm has not terminated after $Q$ iterations})\leq \frac{1}{q!},\]
which converges to 0 as $Q \to \infty$. 

Now we can bound the number of iterations $K$ when the algorithm terminates. Then by we have:
\begin{align*}
    \E[K] &=\sum_{i=0}^\infty \P(K>i)
    \\&\leq\sum_{q=1}^\infty |\ma{Z}|\P(K>|Z|(q-1)+1)+1\\&=\sum_{q=1}^\infty |\ma{Z}|\frac{1}{q!}+1\\&=(e-1)|\ma{Z}|+1
\end{align*}
We have thus  proved that the algorithm terminates in finite number of iterations with probability 1, and furthermore that its expected value is bounded by $(e-1)|\ma{Z}|+1$.

\subsection{Proof of Feasibility and Optimality}

\cite{fletcher1994solving} has shown that   when the cutting plane algorithm terminates, it returns infeasibility or the optimal (feasible) solution. In particular, the cutting plane algorithm returns infeasibility if and only if the original problem is infeasible. Now note that we have assumed the problem is indeed feasible and the formulation of constraints in the stochastic cutting plane algorithm is unchanged from the standard cutting plane algorithm. Therefore, Algorithm \ref{alg:cutting_plane_sto}  returns a feasible solution when it terminates. 

Therefore, we focus on proving that the solution is $\epsilon$-optimal with high probability.

As in the theorem, let $(\bm{z}^*,\bm{\theta}^*)$ be an optimal solution for the original problem and $(\tilde{\bm{z}}^*,\tilde{\bm{\theta}}^*)$ be the solution returned by Algorithm \ref{alg:cutting_plane_sto}. Now let us consider the MILO at the iteration when the algorithm terminated at the $K$th iteration:
\begin{align}
 &\min_{\bm{z}\in \ma{Z}, \bm{\theta}\in \Theta, \eta\geq lb} \;\; \eta\\  \text{s.t.} \quad &\eta \geq f(\bm{z}_i,\bm{\theta}_i; S^n_i)+\nabla f(\bm{z}_i,\bm{\theta}_i; S^n_i)^T{\bm{z}-\bm{z}_i \choose \bm{\theta} - \bm{\theta}_i}, \quad \forall i \in [K],\\ \quad & 0 \geq g_j(\bm{z}_i,\bm{\theta}_i)+\nabla g_j(\bm{z}_i,\bm{\theta}_i)^T{\bm{z}-\bm{z}_i \choose \bm{\theta} - \bm{\theta}_i} ,\quad \forall i \in [K], \; \; \forall j \in  [m] .
\end{align}
For notational simplicity, we define $\eta^*(\bm{z},\bm{\theta})$ to be the optimal value for $\eta$ under the solution $(\bm{z},\bm{\theta})$ for the problem above. In particular, if $(\bm{z},\bm{\theta})$ is feasible, then we have:
\begin{equation}
\eta^*(\bm{z},\bm{\theta}) = \max_{i \in [K]} f(\bm{z}_i,\bm{\theta}_i; S^n_i)+\nabla f(\bm{z}_i,\bm{\theta}_i; S^n_i)^T{\bm{z}-\bm{z}_i \choose \bm{\theta} - \bm{\theta}_i} .
\end{equation}
By construction of the algorithm, there must exist a cutting plane at the optimal point $(\tilde{\bm{z}}^*,\tilde{\bm{\theta}}^*)$. Without loss of generality, we assume it is the final cutting plane, with the form:
\begin{equation}
    \eta \geq f(\tilde{\bm{z}}^*,\tilde{\bm{\theta}}^*; S_K^n)+ \nabla f(\tilde{\bm{z}}^*,\tilde{\bm{\theta}}^*; S^n_K)^T{\bm{z}-\tilde{\bm{z}}^* \choose \bm{\theta} - \tilde{\bm{\theta}}^*}. \label{eq:optimal_const}
\end{equation}
Now we introduce a few lemmas:
\begin{lem}
\label{lem:feasibility}
 $\eta^*(\bm{z}^*,\bm{\theta}^*)\geq f(\tilde{\bm{z}}^*,\tilde{\bm{\theta}}^*; S_K^n)$.
\end{lem}
\proof{Proof:}
By feasibility of  $(\tilde{\bm{z}}^*,\tilde{\bm{\theta}}^*)$, $(\bm{z},\bm{\theta})=(\tilde{\bm{z}}^*,\tilde{\bm{\theta}}^*)$ into Equation (\ref{eq:optimal_const}), the solution returned by  Algorithm \ref{alg:cutting_plane_sto} must satisfy 
\begin{equation}
\eta^*(\tilde{\bm{z}}^*,\tilde{\bm{\theta}}^*)\geq f(\tilde{\bm{z}}^*,\tilde{\bm{\theta}}^*; S^n_K). \label{eq:feasibility}
\end{equation}
Now by optimality of $(\tilde{\bm{z}}^*,\tilde{\bm{\theta}}^*)$, we must have:
\begin{equation}
\eta^*(\bm{z}^*,\bm{\theta}^*)\geq \eta^*(\tilde{\bm{z}}^*,\tilde{\bm{\theta}}^*). \label{eq:optimality}
\end{equation}
Combining equation (\ref{eq:feasibility}) and (\ref{eq:optimality}) gives the required statement.
\QED\endproof
\begin{lem}
\label{lem:lower_bound}
Under the randomization scheme in Algorithm \ref{alg:cutting_plane_sto}, we have, for some absolute constant $D$:
\begin{equation}
    \P\left( \eta^*(\bm{z}^*,\bm{\theta}^*)-f(\bm{z}^*,\bm{\theta}^*;[N])\leq \epsilon \right)\geq  1- 2K\exp\left(-\frac{n\epsilon^2}{(1+\sqrt{p_1+p_2})^2D}\right).
\end{equation}
\end{lem}
\proof{Proof:}
Note that we can write $\eta^*(\bm{z}^*,\bm{\theta}^*)$ as:
\begin{equation}
    \eta^*(\bm{z}^*,\bm{\theta}^*) = \max_{i \in [K]} f(\bm{z}_i,\bm{\theta}_i; S^n_i)+\nabla f(\bm{z}_i,\bm{\theta}_i; S^n_i)^T{\bm{z}^*-\bm{z}_i \choose \bm{\theta}^* - \bm{\theta}_i}. \nonumber 
\end{equation}
Now by Frechet's inequalities:
\begin{align}
    &\P\left( \eta^*(\bm{z}^*,\bm{\theta}^*)-f(\bm{z}^*,\bm{\theta}^*;[N])\leq \epsilon \right)\nonumber \\&= \P\left( \max_{i \in [K]} f(\bm{z}_i,\bm{\theta}_i; S^n_i)+\nabla f(\bm{z}_i,\bm{\theta}_i; S^n_i)^T{\bm{z}^*-\bm{z}_i \choose \bm{\theta}^* - \bm{\theta}_i} -f(\bm{z}^*,\bm{\theta}^*;[N])\leq \epsilon \right)\nonumber \\&=\P\left(  f(\bm{z}_i,\bm{\theta}_i; S^n_i)+\nabla f(\bm{z}_i,\bm{\theta}_i; S^n_i)^T{\bm{z}^*-\bm{z}_i \choose \bm{\theta}^* - \bm{\theta}_i} -f(\bm{z}^*,\bm{\theta}^*;[N])\leq \epsilon \;\;\forall i \in [K]\right)\nonumber \\&\geq \sum_{i=1}^K \P\left( f(\bm{z}_i,\bm{\theta}_i; S^n_i)+\nabla f(\bm{z}_i,\bm{\theta}_i; S^n_i)^T{\bm{z}^*-\bm{z}_i \choose \bm{\theta}^* - \bm{\theta}_i} -f(\bm{z}^*,\bm{\theta}^*;[N])\leq \epsilon \right) - (K-1). \label{eq:full_prob}
\end{align}
We  focus on each term:
\[\P\left( f(\bm{z}_i,\bm{\theta}_i; S^n_i)+\nabla f(\bm{z}_i,\bm{\theta}_i; S^n_i)^T{\bm{z}^*-\bm{z}_i \choose \bm{\theta}^* - \bm{\theta}_i} -f(\bm{z}^*,\bm{\theta}^*;[N])\leq \epsilon \right)\]
By the convexity assumption on $f$ (Assumption \ref{ass:convexity}), we have that:
\begin{equation}
    f(\bm{z}^*,\bm{\theta}^*;[N])\geq f(\bm{z}_i,\bm{\theta}_i; [N])+\nabla f(\bm{z}_i,\bm{\theta}_i; [N])^T{\bm{z}^*-\bm{z}_i \choose \bm{\theta}^* - \bm{\theta}_i}, \quad \forall i \in [K] .\label{eq:convexity}
\end{equation}
Therefore, by substituting Equation (\ref{eq:convexity}), we have, for all $i \in [K]$:
\begin{align}
    &\P\left( f(\bm{z}_i,\bm{\theta}_i; S^n_i)+\nabla f(\bm{z}_i,\bm{\theta}_i; S^n_i)^T{\bm{z}^*-\bm{z}_i \choose \bm{\theta}^* - \bm{\theta}_i} -f(\bm{z}^*,\bm{\theta}^*;[N])\leq \epsilon \right)\\&=  \P\left( (f(\bm{z}_i,\bm{\theta}_i; S^n_i) - f(\bm{z}_i,\bm{\theta}_i; [N]))+\left(\nabla f(\bm{z}_i,\bm{\theta}_i; S^n_i)^T-\nabla f(\bm{z}_i,\bm{\theta}_i; [N])^T\right){\bm{z}^*-\bm{z}_i \choose \bm{\theta}^* - \bm{\theta}_i}\leq  \epsilon \right)\nonumber\\&\geq \P\left( (f(\bm{z}_i,\bm{\theta}_i; S^n_i) - f(\bm{z}_i,\bm{\theta}_i; [N]))\leq \frac{\epsilon}{1+\sqrt{p_1+p_2}} \right. \nonumber \\& \left. \;\;\cap \;\;  \left(\nabla f(\bm{z}_i,\bm{\theta}_i; S^n_i)^T-\nabla f(\bm{z}_i,\bm{\theta}_i; [N])^T\right){\bm{z}^*-\bm{z}_i \choose \bm{\theta}^* - \bm{\theta}_i}\leq  \frac{(\sqrt{p_1+p_2})\epsilon}{1+\sqrt{p_1+p_2}} \right).\label{eq:intermed1}
    \end{align}
    Using the Frechet inequalities, we can further simplify the expression above:
    \begin{align}
    (\ref{eq:intermed1})&\geq\P\left( (f(\bm{z}_i,\bm{\theta}_i; S^n_i) - f(\bm{z}_i,\bm{\theta}_i; [N]))\leq \frac{\epsilon}{1+\sqrt{p_1+p_2}}\right)\\&+ \P\left( \left(\nabla f(\bm{z}_i,\bm{\theta}_i; S^n_i)^T-\nabla f(\bm{z}_i,\bm{\theta}_i; [N])^T\right){\bm{z}^*-\bm{z}_i \choose \bm{\theta}^* - \bm{\theta}_i}\leq  \frac{(\sqrt{p_1+p_2})\epsilon}{1+\sqrt{p_1+p_2}} \right)  -1 .\nonumber &\intertext{Now by Assumption \ref{ass:random}, these expressions lead to :}& \geq 1-\exp\left(-\frac{n\epsilon^2}{(1+\sqrt{p_1+p_2})^2M^2}\right)+ 1-\exp\left(-\frac{n\epsilon^2}{(1+\sqrt{p_1+p_2})^2M'^2}\right)  -1 \\& = 1-2\exp\left(-\frac{n\epsilon^2}{(1+\sqrt{p_1+p_2})^2D}\right),\label{eq:inner_prob}
\end{align}
where $D$ is an absolute constant. Now we substitute Equation (\ref{eq:inner_prob}) into Equation (\ref{eq:full_prob}) to get that:
\begin{align*}
    &\P\left( \eta^*(\bm{z}^*,\bm{\theta}^*)-f(\bm{z}^*,\bm{\theta}^*;[N])\leq \epsilon\right)\\&\geq K- 2K\exp\left(-\frac{n\epsilon^2}{(1+\sqrt{p_1+p_2})^2D}\right)-(K-1)\\&=1- 2K\exp\left(-\frac{n\epsilon^2}{(1+\sqrt{p_1+p_2})^2D}\right).
\end{align*}
as needed. 
\QED\endproof
We now combine the lemmas to derive the final inequality:
\begin{align*}
    &\P\left(f(\tilde{\bm{z}}^*,\tilde{\bm{\theta}}^*;[N])-f(\bm{z}^*,\bm{\theta}^*;[N])\leq \epsilon \right)\\&\geq \P\left(f(\tilde{\bm{z}}^*,\tilde{\bm{\theta}}^*;[N])-f(\tilde{\bm{z}}^*,\tilde{\bm{\theta}}^*;S^n_K)\leq \frac{\epsilon}{2+\sqrt{p_1+p_2}} \; \cap \; f(\tilde{\bm{z}}^*,\tilde{\bm{\theta}}^*;S^n_K)-f(\bm{z}^*,\bm{\theta}^*;N)\leq \frac{(1+\sqrt{p_1+p_2})\epsilon}{2+\sqrt{p_1+p_2}} \right). &\intertext{By the Frechet inequalities, we can simplify the expression above to:}& \geq \P\left(f(\tilde{\bm{z}}^*,\tilde{\bm{\theta}}^*;[N])-f(\tilde{\bm{z}}^*,\tilde{\bm{\theta}}^*;S^n_K)\leq \frac{\epsilon}{2+\sqrt{p_1+p_2}}\right) \\&+ \P\left( f(\tilde{\bm{z}}^*,\tilde{\bm{\theta}}^*;S^n_K)-f(\bm{z}^*,\bm{\theta}^*;N)\leq \frac{(1+\sqrt{p_1+p_2})\epsilon}{2+\sqrt{p_1+p_2}} \right) - 1.&\intertext{Using Lemma \ref{lem:feasibility}, we can bound $f(\tilde{\bm{z}}^*,\tilde{\bm{\theta}}^*;S^n_K)$:}&\geq\P\left(f(\tilde{\bm{z}}^*,\tilde{\bm{\theta}}^*;[N])-f(\tilde{\bm{z}}^*,\tilde{\bm{\theta}}^*;S^n_K)\leq \frac{\epsilon}{2+\sqrt{p_1+p_2}}\right)\\& + \P\left( \eta^*(\bm{z}^*,\bm{\theta}^*)-f(\bm{z}^*,\bm{\theta}^*;N)\leq \frac{(1+\sqrt{p_1+p_2})\epsilon}{2+\sqrt{p_1+p_2}} \right) - 1.&\intertext{Using Assumption \ref{ass:random}, we can evaluate the first expression:}&\geq 1 - \exp\left(\frac{-n\epsilon^2}{(2+\sqrt{p_1+p_2})^2M^2}\right)+ \P\left( \eta^*(\bm{z}^*,\bm{\theta}^*)-f(\bm{z}^*,\bm{\theta}^*;N)\leq \frac{(1+\sqrt{p_1+p_2})\epsilon}{2+\sqrt{p_1+p_2}} \right) - 1.&\intertext{Substituting Lemma \ref{lem:lower_bound} gives us the formula for the second expression:}&\geq 1 - \exp\left(\frac{-n\epsilon^2}{(2+\sqrt{p_1+p_2})^2M^2}\right)+ 1 - 2K\exp\left(\frac{-n\epsilon^2}{(2+\sqrt{p_1+p_2})^2D}\right) - 1\\& = 1- (2K+1)\exp\left(\frac{-n\epsilon^2}{(2+\sqrt{p_1+p_2})^2J}\right),
\end{align*}
where $J$ is an absolute constant, as required. 
\QED
\end{document}